\documentclass[a4paper]{article}

\usepackage[pages=all, color=black, position={current page.south}, placement=bottom, scale=1, opacity=1, vshift=5mm]{background}
\SetBgContents{
	\tt 
}      

\usepackage[margin=1in]{geometry} 

\usepackage{amsmath}
\usepackage{amsthm}
\usepackage{amssymb}

\usepackage[utf8]{inputenc}
\usepackage{hyperref}
\hypersetup{
	unicode,
	pdfauthor={Giuseppe Sirianni},
	pdftitle={phi euler},
	pdfsubject={phi euler},
	pdfkeywords={phi euler},
	pdfproducer={LaTeX},
	pdfcreator={pdflatex}
}  
\usepackage{eso-pic,xcolor}

\usepackage[sort&compress,numbers,square]{natbib}  

\theoremstyle{plain}

\theoremstyle{definition}

\usepackage{graphicx, color}
\graphicspath{{fig/}}

\usepackage{ulem}
\usepackage{contour}
\usepackage{mathtools}
\usepackage{subfig}
\usepackage{float}
\usepackage{authblk}
\usepackage{cancel}
\usepackage{tikz} 

\contourlength{0.8pt}

\renewcommand{\vec}[1]{\boldsymbol{{#1}}}

\newcommand{\normal}{{\boldsymbol{\hat{n}}}}
\newcommand{\divergence}{{\nabla}\cdot}

\newcommand{\partiald}[2]{\dfrac{\partial #1}{\partial #2}}
\newcommand{\partialdd}[2]{\dfrac{\partial^2 #1}{\partial #2^2}}

\newcommand{\veczero}{\vec{0}}
\newcommand{\rhobar}{\overline{\rho}}

\newcommand{\Tbar}{{\overline{T}}}
\newcommand{\varphibar}{{\overline{\varphi}}}
\newcommand{\Ephi}{{E_{\varphi,\rho}}}
\newcommand{\Erho}{{E_{\rho,\varphi}}}
\newcommand{\Ephibar}{{\overline{\Ephi} }}
\newcommand{\Erhobar}{{\overline{\Erho}}}
\newcommand{\Et}{{E^{t}}}

\newcommand*{\var}[1]{\mathord{\mathit{#1}}}

\newcommand{\rhou}{{\var{\rho\vec{u}}}}

\newcommand{\uvector}{{\vec{u}}}

\usepackage{algorithm, algpseudocode} 
\usepackage{mathrsfs} 

\title{An Explicit Primitive Conservative Solver for the Euler Equations with Arbitrary Equation of State}
\author[1,*]{Giuseppe Sirianni}
\author[1]{Alberto Guardone}
\author[1]{Barbara Re}
\author[2]{R\'emi Abgrall} 
\affil[1]{Department of Aerospace Science and Technology, Politecnico di Milano, Via La Masa 34, 20156 Milano, Italy}
\affil[2]{Institute of Mathematics, Universität Zürich, Winterthurerstrasse 190, 8057 Zürich, Switzerland}
\affil[*]{Corresponding author: giuseppe.sirianni@polimi.it}

\date{}
\begin{document}
\maketitle 
\AddToShipoutPictureFG*{
    \AtPageUpperLeft{%
        \raisebox{\dimexpr-\height-2em}{\hspace*{2em}%
        \parbox{1\textwidth}{
                \color{black!30}\bfseries Accepted version of: Computers \& Fluids (2024) 279:106340 \\ 
                Published version available at \href{https://doi.org/10.1016/j.compfluid.2024.106340}{10.1016/j.compfluid.2024.106340} \\
                © 2024. This manuscript version is made available under the \href{https://creativecommons.org/licenses/by-nc-nd/4.0/}{CC-BY-NC-ND 4.0} license 
            }%
        }%
    }%
}%

\section*{Abstract}

This work presents a procedure to solve the Euler equations by explicitly updating, in a conservative manner, a generic thermodynamic variable such as temperature, pressure or entropy instead of the total energy. The presented procedure is valid for any equation of state and spatial discretization. When using complex equations of state such as Span-Wagner, choosing the temperature as the generic thermodynamic variable yields great reductions in the computational costs associated to thermodynamic evaluations. Results computed with a state of the art thermodynamic model are presented, and computational times are analyzed. Particular attention is dedicated to the conservation of total energy, the propagation speed of shock waves and jump conditions. The procedure is thoroughly tested using the Span-Wagner equation of state through the CoolProp thermodynamic library and the Van der Waals equation of state, both in the ideal and non-ideal compressible fluid-dynamics regimes, by comparing it to the standard total energy update and analytical solutions where available.

\subsection*{Keywords}
Primitive formulation, Conservative scheme, Non-ideal compressible fluid dynamics (NICFD), Euler equations, Span-Wagner equation of state

\newpage


\section{Introduction}
The numerical simulation of fluid-dynamics problems is still an extremely active field of research despite thousands of researchers working on the subject for the past century. This is due to a large number of applications, each characterized by its different needs regarding physical modeling, numerical accuracy, and computational costs. The broader subject of fluid-dynamics can be divided into many sub-fields, for example, relativistic~\cite{Rezzolla_2013} and magneto-hydrodynamics (MHD)~\cite{Krause_1980}, incompressible~\cite{Davidson_2021} or compressible~\cite{Balachandran_2006} fluid dynamics, and more. We will concentrate on compressible fluid-dynamics, also named gas-dynamics, which in turn can be subdivided into two further categories, namely ideal and non-ideal~\cite{Guardone_2014,Kluwick_2017,Spinelli_2018} compressible fluid-dynamics (NICFD). Ideal compressible fluid dynamics has been the main interest of the broader CFD research and industry landscape, and it is devoted to the study of flow fields that follow the ideal gas equation of state (EoS)~\cite{Adkins_1983,Cengel_2024}. Starting from the theoretical work of Bethe~\cite{Bethe_1942}, Zel'dovich~\cite{Zeldovich_1946}, and finally Thompson~\cite{Thompson_1971}, a new sub-category of compressible fluid dynamics was born, interested in the study of non-ideal fluids~\cite{Guardone_2014,Kluwick_2017,Spinelli_2018}. Under particular flow conditions, these fluids can exhibit non-classical behaviors such as expansion shocks and compression fans~\cite{Bethe_1942,Zeldovich_1946,Menikoff_1989,Vimercati_2018}. These waves will not occur if the gas is being described through the ideal gas EoS, therefore more accurate EoS such as the Van der Waals (VdW) EoS~\cite{VanDerWaals_1910}, are necessary to describe the thermodynamic behavior of real fluids. Since the first appearance of the VdW EoS, more accurate EoSs have been developed~\cite{Peng_1976,Redlich_1949,Soave_1972}. The current state-of-the-art equation of state is the so-called Span-Wagner type EoS~\cite{Span_1996}, expressed in terms of the Helmholtz free energy. This EoS is at the basis of the thermodynamic library CoolProp~\cite{CoolProp}, which we will use in this work.
The method proposed in this paper will be purposefully developed in an EoS-agnostic manner, that is the proposed strategy is valid for any equation of state. For a comprehensive recap of the latest advancements in the theoretical and numerical modeling of non-ideal compressible fluid dynamics see~\cite{Guardone2024,Gori_2021}.

Non-ideal compressible fluid dynamics is described by the same Euler equations as the ones that describe ideal compressible fluid dynamics. In their conservation form, these are expressed as conservation equations for mass density $\rho$, momentum $\rhou$, and total energy $\Et=E+\rhou^2/2$, with $E$ being the internal energy per unit volume. In NICFD simulations, the most accurate EoS in use is the Span-Wagner type EoS~\cite{Span_1996}, which is expressed as a function of density and temperature $f(\rho, T)$. To evaluate this EoS using two different variables such as density and internal energy $f(\rho,E)$, two different approaches are currently used: look-up tables~\cite{Laughman_2012,Pini_2015,Rubino_2018,Li_2019}, or root searching algorithms~\cite{Gosset_1986,Hickey_2013,Matheis_2018,Trummler_2022}. A similar need also arises in relativistic hydrodynamics where a root-finding algorithm~\cite{Chandra_2017}, or more recently, a machine learning model~\cite{Dieselhorst_2021}, is often employed to retrieve primitive variables from the conservative ones. An alternative and interesting approach that completely replaces the EoS was presented by Saurel et al.~\cite{Saurel_2007} and is based on mechanical relaxation.

In a finite volume CFD code using MUSCL for second-order spatial convergence, one must obtain the thermodynamic state twice for each face separating two control volumes in the mesh. Since in 3D the number of faces scales faster than the number of control volumes, the EoS evaluation costs can add up. In addition to the aforementioned complications introduced by the need to evaluate a more complex EoS, the Riemann solvers used to compute numerical fluxes have also been re-worked to allow for generic EoS~\cite{vinokur_1990,Guardone_2002}.

Analogously to ideal compressible fluid dynamics, conservation~\cite{LeVeque_1992} is a pivotal ingredient in the construction of a robust numerical tool for NICFD problems. Many industrial codes use the finite volume method to discretize the Euler equations due to the conservation properties it naturally embeds. In the compressible regime, the Euler equations are almost universally solved in their conservative form. There exist many different primitive formulations that have been employed in different contexts. For example, Van der Heul et al.~\cite{VanDerHeul_2003} proposed a method that conservatively solves the Euler equations in primitive form with good performance across a large range of Mach numbers. This is done thanks to a particular choice of the pressure scaling and an implicit pressure-correction formula that is only valid for the ideal gas EoS. Toro~\cite{Toro_1998, Toro_2003} described a set of schemes to solve hyperbolic PDEs in non-conservative form and applied them to the Euler equations. Unfortunately, these methods are not conservative by design, and as found by the authors themselves, they tend to compute the wrong shock position. Hughes recapped his works on stabilized methods for compressible fluid-dynamics in~\cite{Hughes_2010}, which also focused on the development of methods using either primitive or entropy variables~\cite{Hauke_1994,Hauke_1998,Hauke_2001}, mostly aimed at the creation of a unified scheme for all-speed flows. As before, these schemes are not conservative by nature. One additional way of dealing with non-conservative formulations is the concept of path-conservation~\cite{Pars_2006,Castro_2008}, although its capability of converging to the exact solution has been questioned~\cite{Abgrall_2010}. Extended Euler systems for multi-phase and multi-component flows have sometimes been presented in primitive form~\cite{Karni_1994,Ogata_1999,Re_2022,Sirianni_2023} as they can be better suited to avoid spurious oscillations across multi-material interfaces~\cite{Abgrall_1996}. A common issue of all primitive formulations is how to preserve conservation. Although they are well-behaved on smooth problems, when discontinuities such as shocks arise in the flow, the loss of conservation can lead to errors in captured shock positions and jump conditions.

Over the years, different approaches have been explored to robustly solve the Euler equations when not using the conservation form. Space-time coupled approaches such as~\cite{Surana_2007,Pesch_2008} look promising but have not been widely adopted. Also, their implementation could require large reworks in well-established CFD codes. Zhang et al.~\cite{Zhang_2018} presented a primitive update of the conservative Euler equations that leverages an implicit dual-time formulation to retrieve exact conservation that was lost by using an approximate linearization of the time derivative to change variables. This approach works well, but it can only recover conservation up to the tolerance on the dual-time integration, which can become expensive for unsteady simulations. Also, Dumbser and Casulli~\cite{Dumbser_2016} presented a semi-implicit conservative solver for generic EoS that updates pressure instead of total energy. This is done by using a Newton-type technique~\cite{Casulli_2013} to solve a system for the unknown pressure on a staggered grid, which could be complicated to implement in already established 2D and 3D unstructured CFD solvers. In the context of shock-free compressible flows, De Michele and Coppola~\cite{Demichele_2023} presented an analysis of different kinetic-energy-preserving schemes, which are currently limited to the ideal gas assumption and the use of central finite difference schemes. More recently, Abgrall~\cite{Abgrall_2023} described a residual correction for the explicit and arbitrarily high-order solution of the primitive Euler equations. This correction works well for any set of variables $[\rho, \vec{u}, \varphi]$, as long as we have a constant derivative of the internal energy with respect to the primitive variable $\varphi$, that is $\Ephi = \left(\partial {E} / \partial {\varphi}\right)_\rho = \textrm{const}$.

In this paper, we will present an update procedure that is explicit in time and allows us to solve the conservative form of the Euler equations by updating $[\rho, \rhou, \varphi]$ instead of $[\rho, \rhou, \Et]$ in a conservative manner without worrying about the spatial discretization. Indeed, we can therefore replace the total energy by any independent generic thermodynamic variable $\varphi$. This procedure does not require the derivation of an evolution equation for $\varphi$ and does not affect any of the properties of its standard total energy update counterpart. The effective choice of $\varphi$ is entirely dependent on the problem at hand. Some applications could require the imposition of boundary conditions on a particular variable such as the pressure ($\varphi=P$), while some others could use a particular EoS where using the enthalpy ($\varphi=h$) is beneficial. We are primarily interested in non-ideal compressible fluid dynamics, so we will showcase results obtained mostly using the temperature ($\varphi=T$) since we can drastically reduce the computational time spent evaluating the thermodynamic properties when using a Span-Wagner type EoS.

The paper is structured as follows: firstly, in section~\ref{sec:spatial}, we will shortly describe the employed spatial discretization. In section~\ref{sec:varphi_update_formula}, we will derive an approximate formula to explicitly update the generic thermodynamic variable $\varphi$ by using the spatial residuals of the total energy equation. Then, in section~\ref{sec:varphibar_secant_search}, we will describe a simple procedure to fix the approximate formula described in the previous section, in order to obtain the exact total energy conservation. To aid the implementation of the presented method, we report a simplified recap in section~\ref{sec:operational_procedure}. Various results are then presented using both the VdW EoS and Span-Wagner type EoS by using the thermodynamic library CoolProp~\cite{CoolProp}. We start with a carbon-dioxide shock tube close to saturation using the Span-Wagner EoS~\cite{Span_1996} in section~\ref{sec:co2_shock} to showcase the method's robustness (section~\ref{sec:co2_shock_comparison}) and the computational costs associated to EoS evaluations (section~\ref{sec:co2_computational_costs}). We then show a dilute nitrogen shock tube with the VdW EoS in section~\ref{sec:n2_shock} to demonstrate how wave speeds and jump conditions are satisfied when using the presented procedure when compared to a simplified approximation. The spatial convergence order is evaluated using the VdW EoS and the Method of Manufactured Solutions (MMS)~\cite{Salari_2000} in section~\ref{sec:mms}. The proposed approach is then tested on multiple dilute gas tests in section~\ref{sec:classical} and dense gas tests in section~\ref{sec:non_classical}  with the VdW EoS, and conclusions are drawn in section~\ref{sec:conclusions}.

\section{Numerical Scheme}\label{sec:numerical_scheme}

The goal of this work is to describe an update procedure that is explicit in time, capable of solving the conservative Euler equations by updating a generic thermodynamic variable $\varphi(x,t)$ instead of the total energy $\Et(x,t)$, in a conservative manner. We will test the method by using $\varphi(x,t)=\left[T, P, e, h, s\right]$, namely temperature, pressure, internal energy, enthalpy, and entropy per unit mass. In short, we will store and update $\varphi$ instead of $\Et$, while still solving the conservative Euler equations. The procedure is structured as follows:
\begin{enumerate}
    \item We compute the spatial discretization of the conservative Euler equations, as described in section~\ref{sec:spatial}
    \item Instead of updating the total energy $\Et$, we use an approximate explicit update formula for $\varphi$, obtained through a linearization which does not conserve total energy, described in section~\ref{sec:varphi_update_formula}
    \item We use a secant root searching algorithm to compute the linearization coefficients present in the aforementioned approximate update formula for $\varphi$, such that total energy is conserved exactly, described in section~\ref{sec:varphibar_secant_search}
\end{enumerate}

As mentioned, the first step to perform this procedure is to compute the spatial residuals of the conservative Euler equations, Eq.~\eqref{eq:euler}, which we can already compute since we are integrating explicitly in time.

\begin{equation}
    \begin{dcases}
        \partiald{\rho}{t}  & + \divergence \left( \rho \vec{u} \right) = 0                   \\
        \partiald{\rhou}{t} & + \divergence \left( \rhou\otimes \vec{u} + P\right) = \veczero \\
        \partiald{\Et}{t}   & + \divergence \left( \left(\Et + P\right)\vec{u} \right) = 0
    \end{dcases}
    \label{eq:euler}
\end{equation}

We will briefly describe the finite volume spatial discretization of the Euler equations we employ in section~\ref{sec:spatial}. We then describe the explicit time update strategy in section~\ref{sec:time}. To simplify the implementation, we recap the procedure in section~\ref{sec:operational_procedure}, forgoing the mathematical details regarding the derivation.

\subsection{Spatial Discretization}\label{sec:spatial}

\begin{figure}[H]
    \centering
    \includegraphics[scale=1]{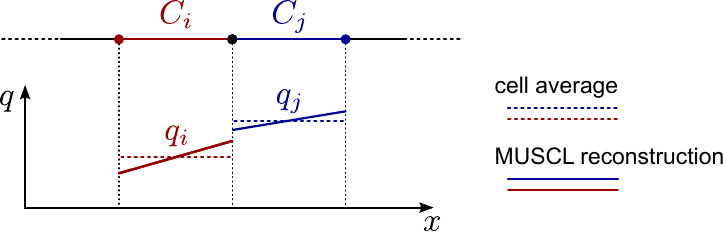}
    \caption{Variable arrangement and MUSCL reconstruction sketch for a generic variable $q$ on two neighbouring control volumes $C_i$ and $C_j$.}
    \label{fig:cv_sketch}
\end{figure}

We will be working in a collocated finite volume solver, with MUSCL reconstruction for second-order spatial convergence (See Fig~\ref{fig:cv_sketch}). We denote with $C_j$ a control volume neighboring control volume $C_i$ and with $A_{i,j}$ their shared area. We also denote with $\sum_{C_j \in \partial {C_i}}$ the sum over all control volumes $C_j$ neighbouring $C_i$.

As mentioned, we will be updating and storing, in each control volume, a generic thermodynamic variable $\varphi(x,t)$ instead of the total energy $\Et(x,t)$. Therefore, the only difference to a standard scheme updating the total energy will be in which variables are stored and reconstructed at each side of the face, namely $\vec{q}_i = \left[\rho_i, \rhou_i, \varphi_i\right]$ and $\vec{q}_j = \left[\rho_j, \rhou_j, \varphi_j\right]$. This means that the numerical fluxes can be computed using any exact or approximate Riemann solver~\cite{Roe_1981,Toro_1994,vinokur_1990,Guardone_2002,Toro_2009}, or also any central difference based schemes~\cite{Jameson_1981}. In this work, we will use an HLLC solver, with the Vinokur \& Montagne~\cite{vinokur_1990} Roe average state for arbitrary EoS. We made this choice because we want to be able to use any equation of state, and also intend to apply the presented scheme to non-ideal compressible fluid dynamics. We name $F^q(\vec{q}_i, \vec{q}_j)$ the numerical flux associated with the conservation equation for an arbitrary variable $q$ from Eq.~\eqref{eq:euler}, on the boundary between two control volume $C_i$ and $C_j$. Using this, we define the spatial residuals $\Phi_i^q$ arising from the discretization of the surface integral of the conservative fluxes for an arbitrary variable $q$ of Eq.~\eqref{eq:euler} on the boundary of a control volume $C_i$ as:
\begin{equation}
    \begin{aligned}
        \Phi_i^{\rho}        & = \sum_{C_j \in \partial {C_i}} A_{i,j} F^{\rho}\left(\vec{q}_i, \vec{q}_j\right)        \\
        \vec{\Phi}_i^{\rhou} & = \sum_{C_j \in \partial {C_i}} A_{i,j} \vec{F}^{\rhou}\left(\vec{q}_i, \vec{q}_j\right) \\
        \Phi_i^{\Et}         & = \sum_{C_j \in \partial {C_i}} A_{i,j} F^{\Et}\left(\vec{q}_i, \vec{q}_j\right)         \\
    \end{aligned}
    \label{eq:residual_cons}
\end{equation}
Note that we compute the residual $\Phi_i^{\Et}$ since we are still interested in solving the total energy conservation equation, although we will update and store the primitive variable $\varphi$. Any high-order spatial reconstruction such as ENO~\cite{Harten_1987} or WENO~\cite{Liu_1994} could be used to reconstruct the left and right states. We will use the MUSCL reconstruction to obtain up to second-order spatial convergence and to avoid spurious oscillations, and we employ the slope limiter by Barth and Jespersen~\cite{Barth_1989}.

\subsection{Time Discretization}\label{sec:time}
In this section, we will describe the conservative time update procedure for $\varphi$. First, we will derive an approximate update formula for $\varphi$ in section~\ref{sec:varphi_update_formula}, obtained through a linearization of the total energy conservation. In section~\ref{sec:varphibar_secant_search} we will then describe a secant root search strategy designed to find the linearization coefficients that recover exact total energy conservation.

As explained in section~\ref{sec:varphi_update_formula}, density $\rho^{n+1}$ and momentum $\rhou^{n+1}$ at the next time step are needed in the approximate update formula for $\varphi$. Since the mass and momentum equations are already in conservative form, we can update them using a forward Euler scheme, that is:
\begin{equation}
    \begin{aligned}
        \rho^{n+1}_i  & = \rho_i^n - \dfrac{\Delta t}{C_i} \Phi^{\rho}_i          \\
        \rhou^{n+1}_i & = \rhou_i^n - \dfrac{\Delta t}{C_i} \vec{\Phi}^{\rhou}_i.
    \end{aligned}
    \label{eq:update_rho_rhou}
\end{equation}
{Note that in this work we are limited to a first order scheme in time, with the same $\var{CFL}$ limitation on the time step size as the corresponding total energy explicit update scheme.

\subsubsection{Approximate update formula for a generic thermodynamic variable \texorpdfstring{$\varphi$}{φ}}\label{sec:varphi_update_formula}

In this section, we will build a formula, similar to what is done in~\cite{Zhang_2018}, which will be used to compute the value of the thermodynamic variable $\varphi_i$ at the next time step, $\varphi_i^{n+1}$, from a given total energy residual $\Phi_i^{\Et}$. This formula does not conserve total energy because it amounts to a linearization of the explicit update~\cite{Zhang_2018}. We will correct this inconsistency through the iterative procedure described in section~\ref{sec:varphibar_secant_search}, while still using an explicit and relatively inexpensive approach.

Denoting the internal energy per unit volume as $E$ and assuming that the equation of state can be expressed in the form $E=E(\rho, \varphi)$, the definition of the total energy reads:
\begin{equation}
    \begin{aligned}
        \Et = E(\rho, \varphi) + \dfrac{1}{2}\rho \vec{u}^2 = \rho e  + \dfrac{1}{2}\rho \vec{u}^2
    \end{aligned}
\end{equation}
where $e$ is the internal energy per unit mass. We can write the derivative of $\Et$ in time as:
\begin{equation}
    \begin{aligned}
        \partiald{\Et}{t} & =\Ephi \partiald{\varphi}{t} +  \Erho \partiald{\rho}{t} + \dfrac{1}{2} \partiald{\rho\vec{u}^2}{t} \\
    \end{aligned}
    \label{eq:t_derivative1}
\end{equation}
where $E_{x,y}$ denotes the partial derivative $\left(\frac{\partial E}{ \partial x}\right)_y$. If we integrate Eq.~\eqref{eq:t_derivative1} in time, between two time steps $t^n$ and $t^{n+1}$ we have:
\begin{equation}
    \begin{aligned}
        \int_{t^n}^{t^{n+1}}\partiald{{\Et}}{{t}}dt & = \int_{t^n}^{t^{n+1}}\left[\Ephi \partiald{\varphi}{t} + \Erho  \partiald{\rho}{t} + \dfrac{1}{2} \partiald{\rho\vec{u}^2}{t}\right]dt                                                      \\[10pt]
        \int_{t^n}^{t^{n+1}}\partiald{{\Et}}{{t}}dt & \simeq \Ephibar \int_{t^n}^{t^{n+1}} \partiald{\varphi}{t} dt + \Erhobar  \int_{t^n}^{t^{n+1}} \partiald{\rho}{t} dt +\dfrac{1}{2} \int_{t^n}^{t^{n+1}}  \partiald{\rho\vec{u}^2}{t} dt  \,.
    \end{aligned}
    \label{eq:discrete_time_jump1}
\end{equation}
In doing this, we made a non-trivial approximation because in principle $\Ephi=\Ephi(\rho(t), \varphi(t))$ and $\Erho = \Erho(\rho(t), \varphi(t))$ are functions of time, whereas we considered them constant and equal to $\Ephibar$ and $\Erhobar$, respectively. This linearization is why the update formula described in this section is not exactly conservative. In section~\ref{sec:varphibar_secant_search} we will outline a way of selecting their values to reconcile the approximate and exact jumps of internal energy.

Introducing the notation $\Delta x = x^{n+1} - x^n$ for the jump of a variable $x$ during a time step, we can write Eq.~\eqref{eq:discrete_time_jump1} in a compact way as:
\begin{equation}
    \begin{aligned}
        \Delta {\Et} & \simeq \Ephibar  \Delta{\varphi} + \Erhobar    \Delta {\rho} + \dfrac{1}{2}\Delta {\rho}{\vec{u}}^2.
    \end{aligned}
    \label{eq:discrete_time_jump2}
\end{equation}

Let us notice that in Eq.~\eqref{eq:discrete_time_jump2} there are  two distinct contributions to the jump of total energy, namely the approximation of the internal energy jump $\Delta E^{\textrm{approx}}$ and the exact kinetic energy jump $\Delta \rhou^2 / 2$ that we already know because we have already updated the density and momentum through Eq.~\eqref{eq:update_rho_rhou}:
\begin{equation}
    \begin{aligned}
        \Delta {\Et} & \simeq \underbrace{\Ephibar  \Delta{\varphi} + \Erhobar    \Delta {\rho}}_{\Delta E^{\textrm{approx}}} + \underbrace{\dfrac{1}{2}\Delta {\rho}{\vec{u}}^2	}_{\text{exact and known}} \!\!\!\!\!\!\!.
    \end{aligned}
    \label{eq:delta_E_approx1}
\end{equation}

We can now use the approximate jump of total energy from Eq.~\eqref{eq:discrete_time_jump2} to write an approximation of the conservation of discrete total energy for control volume $C_i$:
\begin{equation}
    \begin{aligned}
        \int_{{C_i}} \Delta {\Et_i}  d{V}                                                                                                  & + \Delta {t} {\Phi}_{i}^{{\Et}}  = 0    \\
        \int_{{C_i}} \left(\Ephibar_i\Delta \varphi_i + \Erhobar_i \Delta {\rho}_i  + \dfrac{1}{2}   \Delta{\rho}{\vec{u}}^2_i \right)d{V} & + \Delta {t} {\Phi}_{i}^{{\Et}} = 0 \,.
    \end{aligned}
    \label{eq:cons_tot_energy1}
\end{equation}

For the sake of readability, we will drop subscript $(\cdot)_i$ denoting the control volume because we are only looking at a single control volume $C_i$. Let us split the approximate total energy conservation from Eq.~\eqref{eq:cons_tot_energy1} into two parts, and rename $B$ the one containing $\Delta \varphi$:
\begin{equation}
    \underbrace{ \int_{{C}}  \Ephibar \Delta \varphi d{V}}_{B}+  \int_{{C}} \left( \Erhobar \Delta {\rho}  + \dfrac{1}{2}   \Delta{\rho}{\vec{u}}^2 \right)d{V} +  \Delta {t} {\Phi}_{i}^{{\Et}}  = 0 \,.
    \label{eq:cons_tot_energy2}
\end{equation}

We will use the same forward Euler time update scheme for the generic thermodynamic variable $\varphi$ as the one used for density and momentum in Eq.~\eqref{eq:update_rho_rhou}. We therefore define a pseudo-spatial residual $\Phi^\varphi$ in control volume $C$ that satisfies exactly the approximate total energy conservation, Eq.~\eqref{eq:cons_tot_energy2}
\begin{equation}
    \varphi^{n+1} = \varphi^n - \dfrac{\Delta t}{C} \Phi^\varphi\,,
    \label{eq:delta_phi1}
\end{equation}
which we can easily invert to obtain:
\begin{equation}
    \Delta \varphi = - \dfrac{\Delta t}{C} \Phi^\varphi \,.
    \label{eq:delta_phi2}
\end{equation}

Substituting Eq.~\eqref{eq:delta_phi2} into $B$ from Eq.~\eqref{eq:cons_tot_energy2}, and introducing $\overline{\omega}$, we obtain
\begin{equation}
    \begin{aligned}
        B & = -\Delta t \Phi^\varphi  \dfrac{1}{C}\int_{{C}}\Ephibar  d{V} \\
          & = -\Delta t \Phi^\varphi  \overline{\omega}   \,.              \\
    \end{aligned}
    \label{eq:b_varphi}
\end{equation}
Let us now substitute this expression for $B$ back into the approximate total energy conservation Eq.~\eqref{eq:cons_tot_energy2}:
\begin{equation}
    -    \Delta t \Phi^\varphi \overline{\omega}    +  \int_{{C}} \left( \Erhobar \Delta {\rho}  + \dfrac{1}{2}   \Delta{\rho}{\vec{u}}^2 \right)d{V} +  \Delta {t} {\Phi}^{{\Et}}   = 0 \,.
    \label{eq:cons_tot_energy3}
\end{equation}
from which we can explicitly express the pseudo-spatial residual $\Phi^\varphi$ for the generic thermodynamic variable $\varphi$:
\begin{equation}
    \Phi^\varphi = \dfrac{1}{\overline{\omega}  \Delta t } \left[ \int_{{C}} \left( \Erhobar \Delta {\rho}  + \dfrac{1}{2}   \Delta{\rho}{\vec{u}}^2 \right)d{V} +  \Delta {t} {\Phi}^{{\Et}}\right] \,.
    \label{eq:varphi_residual}
\end{equation}
Therefore the final form of the explicit update formula for a generic thermodynamic variable $\varphi$  reads
\begin{equation}
    \begin{dcases}
        \Phi^\varphi = \dfrac{1}{\overline{\omega}  \Delta t } \left[ \int_{{C}} \left( \Erhobar \Delta {\rho}  + \dfrac{1}{2}   \Delta{\rho}{\vec{u}}^2 \right)d{V} +  \Delta {t} {\Phi}^{{\Et}}\right] \\[10pt]
        \varphi^{n+1} = \varphi^n - \dfrac{\Delta t}{C} \Phi^\varphi \,.
    \end{dcases}
    \label{eq:varphi_update_formula_final}
\end{equation}

The presence of $\Delta \rho = \rho^{n+1}-\rho^n$ and $\Delta \rhou^2 = \left(\rhou^2\right)^{n+1}-\left(\rhou^2\right)^n$ in Eq.~\eqref{eq:varphi_update_formula_final} highlights the need to perform the update of $\rho$ and $\rhou$ before updating $\varphi$. Notice that in Eq.~\eqref{eq:varphi_update_formula_final} the pseudo-spatial residual $\Phi^\varphi$ is a function of unknowns $\Ephibar$ and $\Erhobar$. As a consequence, also the updated value of the thermodynamic variable $\varphi^{n+1}$ will be a function of these unknowns:
\begin{equation}
    \begin{dcases}
        \Phi^\varphi                      = \Phi^\varphi(\Ephibar, \Erhobar)                                 \\
        \varphi^{n+1}(\Ephibar, \Erhobar) = \varphi^n - \dfrac{\Delta t}{C}\Phi^\varphi(\Ephibar, \Erhobar). \\
    \end{dcases}
    \label{eq:phi_update_function_bar}
\end{equation}
Furthermore, also $\Delta E^{\textrm{approx}}$,  as defined in Eq.~\eqref{eq:delta_E_approx1}, depends on $\Ephibar$ and $\Erhobar$.

Note that the definition of local conservation remains the same regardless of the time discretization employed, since the time jumps $\Delta y \equiv \int_{t^n}^{t^{n+1}} y \var{dt}$ that appear in Eq.~\eqref{eq:discrete_time_jump1} are exact. One could apply this notion to extend the presented time scheme to higher orders.

\subsubsection{How to retrieve total energy conservation -- \texorpdfstring{$\varphibar$}{φ} search}\label{sec:varphibar_secant_search}

In section~\ref{sec:varphi_update_formula} we derived a formula that allows us to update a generic thermodynamic variable $\varphi$ instead of the total energy. The derived formula in Eq.~\eqref{eq:varphi_update_formula_final} is not exactly conservative due to the approximation made in Eq~\eqref{eq:discrete_time_jump2}, where we introduced the constant approximate values $\Ephibar$ and $\Erhobar$. The goal of this section is to show a simple procedure that can be used to select the values of $\Ephibar$ and $\Erhobar$, such that they enforce exact total energy conservation. Let us first recall the formula for the approximate internal energy jump $\Delta E^{\textrm{approx}}$ from Eq.~\eqref{eq:delta_E_approx1} that has been used to derive the $\varphi$ update formula Eq.~\eqref{eq:varphi_update_formula_final}. We can notice that it is a function of the  constant approximate values $\Ephibar$ and $\Erhobar$, which are unknown:
\begin{equation}
    \Delta E^{\textrm{approx}}(\Ephibar, \Erhobar) = \Ephibar  \Delta{\varphi} + \Erhobar    \Delta {\rho} \,.
    \label{eq:delta_E_approx2}
\end{equation}

We begin by imposing thermodynamic consistency at an unknown approximate thermodynamic state $\left(\rhobar,\varphibar\right)$:
\begin{equation}
    \begin{dcases}
        \Ephibar & = \Ephi (\rhobar, \varphibar)     \\
        \Erhobar & = \Erho (\rhobar, \varphibar) \,.
    \end{dcases}
    \label{eq:thermo_consistency}
\end{equation}
In doing this, we have changed the unknowns of our problem from $\Ephibar$ and $\Erhobar$ to $\rhobar$ and $\varphibar$. Given the fact that we already know the exact updated density and that its change between two time steps is linear in $\Delta t$, we assume that we can choose the value of $\rhobar$ as
\begin{equation}
    \rhobar = \dfrac{1}{2}\left(\rho^{n+1} + \rho^n\right) \,.
    \label{eq:rhobar}
\end{equation}
Since now $\rhobar$ is fixed, the problem is reduced to a single unknown $\varphibar$. Recalling Eq.~\eqref{eq:phi_update_function_bar}, we can see that by using the $\varphi$ update formula Eq.~\eqref{eq:varphi_update_formula_final}, the updated thermodynamic variable $\varphi^{n+1}$ at the next time step is now also a function of the unknown $\varphibar$:
\begin{equation}
    \varphi^{n+1}(\varphibar) \;\;\;\;\;\;\; \rightarrow \;\;\;\;\;\;\; \Delta\varphi(\varphibar)= \varphi^{n+1}(\varphibar) - \varphi^n \,.
\end{equation}

If we now use $\varphi^{n+1}(\varphibar)$ to write both the approximate internal energy jump $\Delta E^{\textrm{approx}}$ from Eq.~\eqref{eq:delta_E_approx1} and its thermodynamically exact counterpart $\Delta E^{\textrm{tmd}}$ through the EoS we obtain
\begin{equation}
    \begin{dcases}
        \Delta E^{\textrm{approx}}(\varphibar) & = \Ephi(\rhobar,\varphibar)  \Delta\varphi(\varphibar) + \Erho(\rhobar,\varphibar)    \Delta {\rho} \\[10pt]
        \Delta E^{\textrm{tmd}}(\varphibar)    & = E(\rho^{n+1}, \varphi^{n+1}(\varphibar)) - E(\rho^n, \varphi^n) \,.
    \end{dcases}
    \label{eq:delta_E_approx3}
\end{equation}
Ideally, we would want these two internal energy jumps to coincide. We therefore want to find $\varphibar$ such that $\Delta E^{\textrm{tmd}} (\varphibar) = \Delta E^{\textrm{approx}} (\varphibar)$. This would entail exact total energy conservation, since the jump in kinetic energy $ \Delta \rhou^2 / 2$ is already known exactly. We can construct a scalar function $F(\varphibar)$ such that its root corresponds to the value of $\varphibar$ we are looking for:
\begin{equation}
    F(\varphibar) = \dfrac{\Delta E^{\textrm{tmd}}(\varphibar) - \Delta E^{\textrm{approx}}(\varphibar)}{E(\rho^n, \varphi^n)}
    \label{eq:F_definition}
\end{equation}
where the division by $E(\rho^n, \varphi^n)>0$ is performed to scale $F(\varphibar)$. The problem is now to find a root of $F(\varphibar)$ for which we will use the secant root search method, where we iterate on $\varphibar_{(k)}$ until:
\begin{equation}
    \lvert F(\varphibar_{(k)}) \rvert < \var{tol} = 10^{-14} \,.
    \label{eq:F_secant}
\end{equation}
We chose $10^{-14}$ as it is close enough to machine precision without being too restrictive. Looser tolerances will yield slightly lower iteration counts and cause a loss in conservation, since $F$ is an error metric directly related to the total energy conservation error. There is no clear choice for the bounds of the secant search. Given the fact that the size of the time-step is already bound by the CFL condition, in this work we use $\left[\varphi^n, \varphi^{n+1}\right]$. Since, of course, $\varphi^{n+1}$ is unknown, we use a first guess for the updated generic thermodynamic variable $\tilde{\varphi}^{n+1}$ computed evaluating Eq.~\ref{eq:varphi_update_formula_final} with $\rhobar=(\rho^{n+1}+\rho^n)/2$ and $\varphibar=\varphi^n$.

We have to point out some possible theoretical shortcomings in the proposed approach. To the best of our knowledge, there is no clear way to prove that such a root exists for a generic EoS and that if it does, it is the only root. In our practical experience, $F(\varphibar)$ is always linear with respect to $\varphibar$ for any choice of $\varphi=\left[T, P, e, h, s, \dots\right]$, see Fig.~\ref{fig:n2_shock_F}, therefore the secant search algorithm always converges in at most seven iterations, but usually one or the initial guess are enough. When the change in the internal energy is close to zero, there can be a situation where there is no root but with
\begin{equation}
    \lvert F(\varphibar) \rvert  \sim \epsilon^{\textit{machine}}\;\;\;\;\;\forall \varphi \in \left[\varphi^n,\tilde{\varphi}^{n+1}\right] \,.
    \label{eq:F_machine_precision}
\end{equation}
In this situation, we can consider conservation as satisfactorily met by the initial guess, since it is close to machine precision. We were unable to find test cases where this causes any appreciable loss in total energy conservation.

\begin{figure}
    \centering
    \subfloat[Temperature, $\varphi=T$]{
        \includegraphics[scale=1]{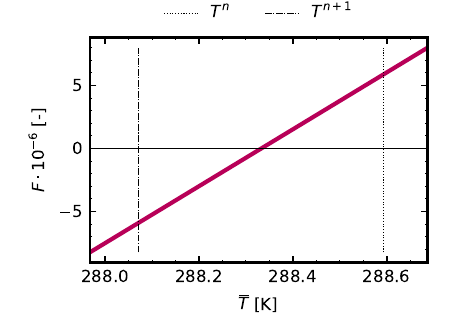}
        \label{subfig:n2_shock_F_temperature}
    }%
    \subfloat[Pressure, $\varphi=P$]{
        \includegraphics[scale=1]{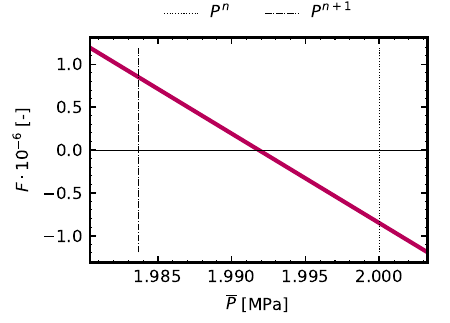}
        \label{subfig:n2_shock_F_pressure}
    }\\
    \subfloat[Specific energy, $\varphi=e$]{
        \includegraphics[scale=1]{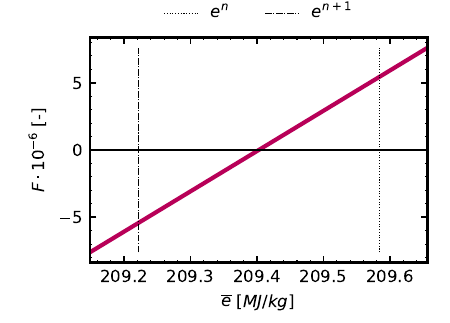}
        \label{subfig:n2_shock_F_energy}
    }%
    \subfloat[Specific enthalpy, $\varphi=h$]{
        \includegraphics[scale=1]{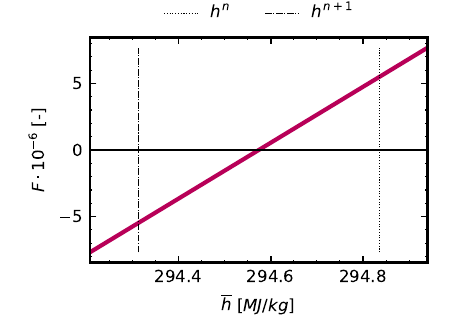}
        \label{subfig:n2_shock_F_enthalpy}
    }\\
    \subfloat[Specific entropy, $\varphi=s$]{
        \includegraphics[scale=1]{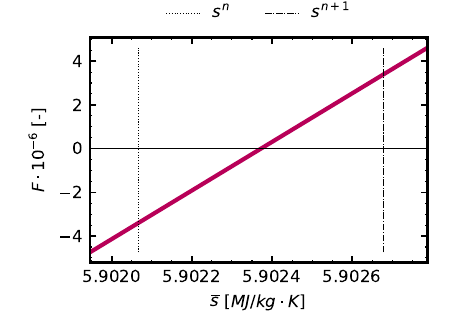}
        \label{subfig:n2_shock_F_entropy}
    }
    \caption{Example plots of $F(\varphibar)$ for various $\varphi$ choices at the initial time step discontinuity of a Nitrogen shock tube using the Span-Wagner EoS~\cite{Span_2000}.}
    \label{fig:n2_shock_F}
\end{figure}

\subsubsection{Operational Recap of the Method}\label{sec:operational_procedure}

In this section, we quickly recap the conservative update procedure for $\varphi$, to simplify its implementation. First, for the sake of clarity, let us rewrite Eq.~\eqref{eq:varphi_update_formula_final} by explicitly stating the thermodynamic dependencies on $\rhobar$ and $\varphibar$:
\begin{equation}
    \begin{dcases}
        \Phi^\varphi = \dfrac{1}{\overline{\omega} (\Ephi(\rhobar, \varphibar)) \Delta t } \left[ \int_{{C}} \left( \Erho(\rhobar, \varphibar) \Delta {\rho}  + \dfrac{1}{2}   \Delta{\rho}{\vec{u}}^2 \right)d{V} +  \Delta {t} {\Phi}^{{\Et}}\right] \\[10pt]
        \varphi^{n+1} = \varphi^n - \dfrac{\Delta t}{C} \Phi^\varphi
    \end{dcases}
    \label{eq:varphi_update_formula_operational}
\end{equation}

We also recall the secant formula:
\begin{equation}
    \varphibar_{(k+1)} = \varphibar^{-}_{(k)} - F^{-}_{(k)} \dfrac{\varphibar^{+}_{(k)} - \varphibar^{-}_{(k)}}{F^{+}_{(k)} - F^{-}_{(k)}}
    \label{eq:secant}
\end{equation}

The main steps of the conservative update procedure for $\varphi$ are the following.
\begin{enumerate}
    \item Using the solution $\left[\rho^n, \rhou^n, \varphi^n\right]$ at the previous time step $t^n$, compute the spatial residuals of the conservative Euler equations Eq.~\eqref{eq:euler}, namely $\left[\Phi^{\rho},\vec{\Phi}^{\rhou},\Phi^{\Et}\right]$ as defined in Eq.~\eqref{eq:residual_cons}, using the preferred choice of numerical fluxes or central difference scheme.
    \item Update the density and momentum to compute $\rho^{n+1}$ and $\rhou^{n+1}$ at the next time step using Eq.~\eqref{eq:update_rho_rhou}, as the mass and momentum equations are already in conservative form.
    \item Compute an initial guess $\tilde{\varphi}^{n+1}$ for $\varphi^{n+1}$ at the next time step $t^{n+1}$. In this work we use Eq.~\eqref{eq:varphi_update_formula_operational}, with $\rhobar=(\rho^{n+1} + \rho^{n})/2$ and $\varphibar=\varphi^n$.
    \item According to the definition of $F(\varphibar)$ from Eq.~\eqref{eq:F_definition}, compute the initial bounds for the secant search $F^{+}_{(0)} = \max\left[F(\varphi^n), F(\tilde{\varphi}^{n+1})\right]$ and $F^{-}_{(0)} = \min\left[F(\varphi^n), F(\tilde{\varphi}^{n+1})\right]$, and also save the maximum and minimum bounds  $\varphibar^{+}_{(0)}$ and $\varphibar^{-}_{(0)}$ accordingly. Other choices for the initial guess $\tilde{\varphi}^{n+1}$ are possible.
    \item Begin an iterative loop, with a user-defined tolerance and maximum number of iterations.
          \begin{enumerate}
              \item Compute $\varphibar_{(k+1)}$ at iteration $k+1$ using the secant formula, see Eq.~\eqref{eq:secant}.
              \item Using $\varphibar_{(k+1)}$ and Eq.~\eqref{eq:varphi_update_formula_operational} compute the next guess value for $\varphi^{n+1}_{(k+1)}$.
              \item Using $\varphibar_{(k+1)}$ and $\varphi^{n+1}_{(k+1)}$, compute $F(\varphibar_{(k+1)})$ and check for convergence, namely $\lvert F(\varphibar_{(k+1)}) \rvert < \textit{tol}$.
              \item If convergence has not been reached, save the new values for the secant search $F^{+}_{(k+1)}$, $F^{-}_{(k+1)}$, $\varphibar^{+}_{(k+1)}$ and $\varphibar^{-}_{(k+1)}$ and perform another iteration.
          \end{enumerate}
    \item Save $\varphi^{n+1}_{(k+1)}$ as the next time step value of $\varphi$.
\end{enumerate}

The scheme described here has been derived without assuming a 1D setting.
Hence, its extension to 2D or 3D schemes does not pose significant difficulties. In particular, we only need to compute the jump in kinetic energy, $\Delta (\rho \uvector\cdot\uvector)/2$, and the spatial residual of the total energy equation, $\Phi^{\Et}$, using the multi-dimensional discretization of choice. Then, the proposed time update scheme can be applied as is.

\section{Results}

In this section, we will showcase a variety of test cases. First, we use the Span-Wagner EoS~\cite{Span_1996} through the thermodynamic library CoolProp~\cite{CoolProp} in a carbon-dioxide shock tube test close to saturation in section~\ref{sec:co2_shock}. Here, we compare solution quality (section~\ref{sec:co2_shock_comparison}) and EoS computational costs (section~\ref{sec:co2_computational_costs}) for various choices of the generic thermodynamic variable $\varphi=\left[T,P,e,h,s\right]$. We then use the VdW EoS in a dilute nitrogen shock tube in section~\ref{sec:n2_shock} to show how the proposed procedure is necessary to retrieve exact total energy conservation, and how failing to do so affects captured shock speeds and the fulfillment of jump conditions. In section~\ref{sec:mms} we perform order testing using the VdW EoS and the method of manufactured solutions (MMS)~\cite{Roache_1998}. Following this, we show a series of tests using the VdW EoS in the ideal regime in section~\ref{sec:classical} and finally in the NICFD regime in section~\ref{sec:non_classical}.

The VdW EoS has been chosen because it is capable of describing non-ideal behavior. See appendix~\ref{appendix:vdw} for more details on the VdW EoS. We compare results to exact solutions of the Riemann problems computed following the approach described in Quartapelle et al.~\cite{Quartapelle_2003}. For reproducibility reasons, we need to define here the constants $\delta = \gamma - 1$ where $\gamma$ is the ratio of specific heats, $a$ which describes the magnitude of the attractive forces between molecules, $b$ which is the co-volume and $R$ which is the gas constant.

Since one of the main goals of this work is assessing conservativity, we will often show plots of the dimensionless imbalance of total energy, named $I^{\Et}$, which we define as
\begin{equation}
    I^{\Et}(t) = \dfrac{1}{C_{\Omega}\left[\int_{\Omega} \Et(\tau) dV\right]_{\tau=0} }\left\lbrace\left[\int_{\Omega} \Et(\tau) dV\right]_{\tau=t} \!\!\! -\left[\int_{\Omega} \Et(\tau) dV\right]_{\tau=0} \!\!\! +\int_{\tau=0}^{\tau=t}\left[\oint_{\partial \Omega}\!\!\! \vec{F}^{\Et}(\tau)\cdot\normal dA\right]d\tau\right\rbrace
\end{equation}
The $(\cdot)^t$ superscript in $\Et$ does not refer to time $t$. $I^{\Et}$ is both normalized in terms of the domain's volume $C_{\Omega}$ and of the integral of the initial total energy $\int_{\Omega}\Et(0) dV$, and is therefore dimensionless. Also, it is a time-varying quantity as it measures the amount of total energy generated or destroyed by the method in the domain since $t=0$.

We also compare the results of the proposed scheme against numerical simulations from a standard finite volume solver for the Euler equations that is identical to the previous one (same numerical fluxes, limiting, compiler) except that, instead of updating the thermodynamic variable $\varphi$ using the method presented in this paper, we update the total energy $\Et$ explicitly as
\begin{equation}
    {\Et}_i^{n+1} = {\Et}_i^n - \dfrac{\Delta t}{C_i} \Phi^{\Et}_i .
\end{equation}
We use {``${\varphi\text{ Update}}$"} to denote the results obtained with the update of the thermodynamic variable $\varphi$ presented in this work. In contrast, we use {``${\Et\text{ Update}}$"} to denote the standard finite volume solver. For the latter, we also compare results obtained using both the analytical EoS ({``${\Et\text{ Update}}$"}) or the computationally efficient look-up tables (LUT) as implemented in CoolProp ({``${E^{t,\text{(LUT)}}\text{ Update}}$"}).

Presenting the thermodynamic results, we will also refer to the fundamental derivative of gas-dynamics $\Gamma$~\cite{Thompson_1971} and the compressibility factor $Z$, which are defined as:
\begin{equation}
    \Gamma(s,v) = \dfrac{v^3}{2 c^2} \partialdd{P(s,v)}{v} ,\;\;\;\;\;\;\;\;\;\;\;\;\;\;\;\;\;\; Z = \dfrac{P}{\rho R T} \,.
    \label{eq:G_definition}
\end{equation}

\subsection{Carbon-Dioxide Shock Tube (Span-Wagner)}\label{sec:co2_shock}

In this section we show the results for a carbon-dioxide shock tube, see table~\ref{tab:co2_shock_data} for the data and Fig.~\ref{subfig:co2_shock_saturation} for the thermodynamic path on the $(P,\rho)$ plane. We use this test to assess the behavior and performance of the method in a worst case scenario, with thermodynamic states spanning from supercritical to subcritical conditions, while also approaching the saturation line. The solution is composed of a left-running rarefaction fan that expands from a supercritical state up to the vapour saturation curve, and a right running compression shock. In particular, we show the difference between different possible choices for $\varphi$ using the Span-Wagner EoS~\cite{Span_1996} through the thermodynamic library CoolProp~\cite{CoolProp}. We devised this test to stay close to the saturation curve (see Fig.~\ref{subfig:co2_shock_saturation}) to show that the proposed method works for any variable, even when far from the dilute gas regime (see the compressibility factor in Fig.~\ref{subfig:co2_shock_Z}). See appendix~\ref{appendix:n2_shock} for an analogous analysis for a dilute nitrogen shock tube test also using CoolProp.

\begin{table}[H]
    \renewcommand{\arraystretch}{1.2}
    \centering
    \begin{tabular}{l | l | l l}
        \textit{Equation of State}   & \textit{Numerical}                   & ${x < 0.5\; \mathrm{m}}$     & ${x > 0.5\; \mathrm{m}}$     \\ \hline
        Span-Wagner~\cite{Span_1996} & $x\in\left[0,1\right]\; \mathrm{m}$  & $\rho=350 \;\mathrm{kg/m^3}$ & $\rho=100 \;\mathrm{kg/m^3}$ \\
        Carbon-Dioxide               & $t_{\var{final}}=0.001 \;\mathrm{s}$ & $u = 0 \; \mathrm{m/s}$      & $u = 0 \; \mathrm{m/s}$      \\
                                     & $\var{CFL}=0.1$                      & $P = 12 \; \mathrm{{MPa}}$   & $P = 4 \; \mathrm{{MPa}}$
    \end{tabular}
    \caption{Carbon-dioxide shock tube test initial conditions, thermodynamic, domain and numerical data.}
    \label{tab:co2_shock_data}
\end{table}

\begin{figure}[H]
    \centering
    \subfloat[Compressibility Factor]{
        \includegraphics[scale=1]{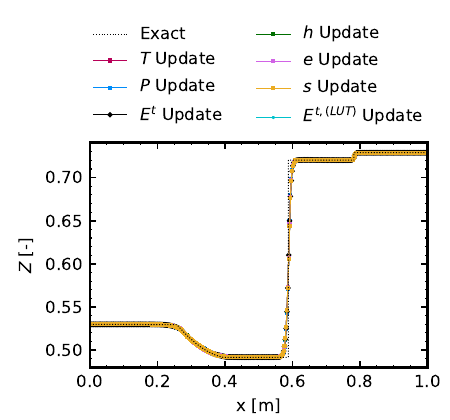}
        \label{subfig:co2_shock_Z}
    }%
    \subfloat[Thermodynamic Plane]{
        \includegraphics[scale=1]{
            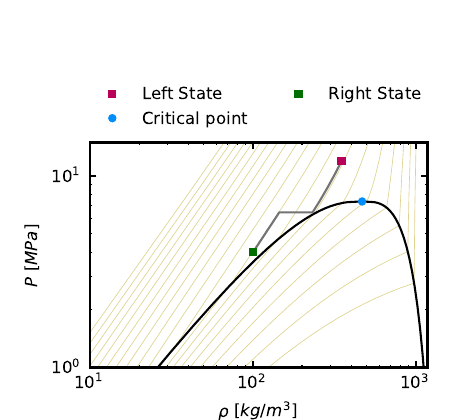}
        \label{subfig:co2_shock_saturation}
    }
    \caption{Carbon-dioxide shock tube test with $\Delta x = 0.025 \;\mathrm{m}$ using the Span-Wagner EoS~\cite{Span_1996} through the CoolProp thermodynamic library~\cite{CoolProp}. Compressibility factor profile and numerical solution on the on the pressure-density thermodynamic plane at $t=0.001\;\mathrm{s}$, with isentropes in yellow and the vapour saturation curve in black.}
    \label{fig:co2_shock_thermo}
\end{figure}

\subsubsection{Comparison of \texorpdfstring{$\varphi$}{φ} choices}\label{sec:co2_shock_comparison}

In Fig.~\ref{fig:co2_shock_comp} we show profiles of density, pressure, and velocity compared to the exact solution. There are negligible variations in the results between all choices of $\varphi$, with the only exception being the specific entropy $s$ in yellow, showing a small kink on the rightmost part of the rarefaction wave. This is probably due to the MUSCL reconstruction of $s$ near the saturation vapor curve: if we look at first-order results without MUSCL reconstruction in Fig.~\ref{fig:co2_shock_comp_1st_order}, we can qualitatively see that all $\varphi$ choices perfectly overlap. To confirm this quantitatively, we report the dimensionless difference of density and pressure profiles with respect to the results of the standard $\Et$ update in Fig.~\ref{subfig:co2_shock_err_1st_p} and~\ref{subfig:co2_shock_err_1st_rho}. Here we can also appreciate how the proposed approach yields results that are much closer to the standard $\Et$ update than using look-up tables. All other results with MUSCL in Fig.~\ref{fig:co2_shock_comp} are in great agreement with the solution obtained using the standard $\Et$ update. More importantly, if we focus on the total energy imbalance $I^{\Et}$ in Fig.~\ref{subfig:co2_shock_comp_energy_conservation_log} we see how all conserve total energy close to machine precision. A closer look reveals slight differences in the total energy imbalances of different choices of $\varphi$. Namely, the specific internal energy $e$ behaves identically to the standard $\Et$ update (both using the EoS directly and using look-up-tables). Temperature $T$, specific enthalpy $h$, and specific entropy $s$ show some similarities, while pressure $P$ showcases some sporadic steps that may at first glance appear large. If we assume a loss of conservation per element close to machine precision (for our setup $\epsilon^{\textit{machine}}\simeq 2\cdot10^{-16}$), on a $400$ element mesh the total loss per time step would be approximately $400 \cdot \left(2 \cdot 10^{-16} \right)\simeq 8 \cdot 10^{-14}$. The simulation ran for $2920$ iterations, therefore we can estimate an approximate loss of conservation due to truncation errors accumulating for the total simulation as $2920 \cdot 8 \cdot 10^{-14} \simeq 2.4 \cdot 10^{-10}$, which is very comparable to the final total energy imbalance for the $P$ update of $\sim 2\cdot 10^{-11}$. This estimate is a worst-case scenario, as errors could cancel out in the domain instead of accumulating.

\begin{figure}
    \centering
    \subfloat[Density]{
        \includegraphics[scale=1]{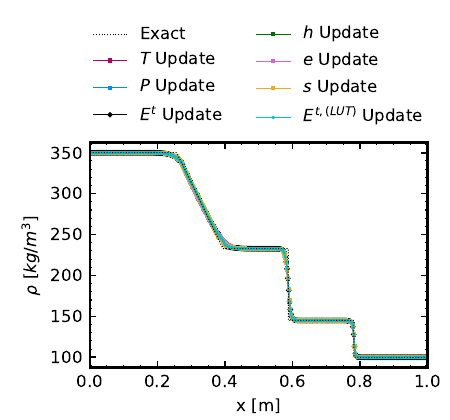}
        \label{subfig:co2_shock_comp_rho}
    }%
    \subfloat[Pressure]{
        \includegraphics[scale=1]{
            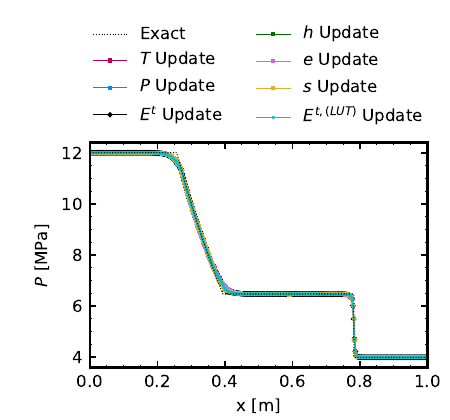}
        \label{subfig:co2_shock_comp_p}
    }\\
    \subfloat[Velocity]{
        \includegraphics[scale=1]{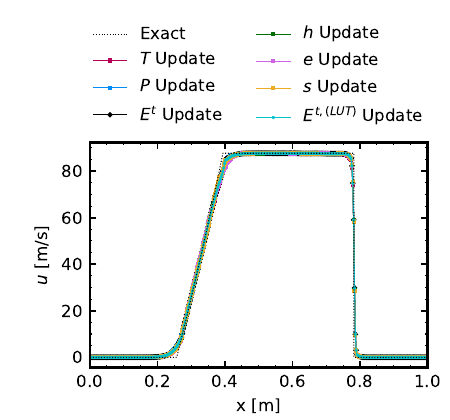}
        \label{subfig:co2_shock_comp_u}
    }%
    \subfloat[Total energy imbalance]{
        \includegraphics[scale=1]{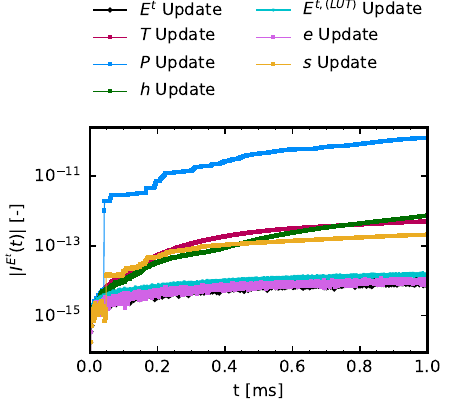}
        \label{subfig:co2_shock_comp_energy_conservation_log}
    }
    \caption{Carbon-dioxide shock tube test results using MUSCL with $\Delta x = 0.025\;\mathrm{m}$ using the Span-Wagner EoS~\cite{Span_1996} through the CoolProp thermodynamic library~\cite{CoolProp}. Density, pressure, and velocity profiles at $t=0.001\;\mathrm{s}$ and absolute value of the total energy imbalance in time. Comparison between the $\varphi$ update scheme for temperature, pressure, specific energy, specific enthalpy, specific entropy, and the standard total energy update (both using the analytical EoS and a look-up table).}
    \label{fig:co2_shock_comp}
\end{figure}

\begin{figure}
    \centering
    \subfloat[Density]{
        \includegraphics[scale=1]{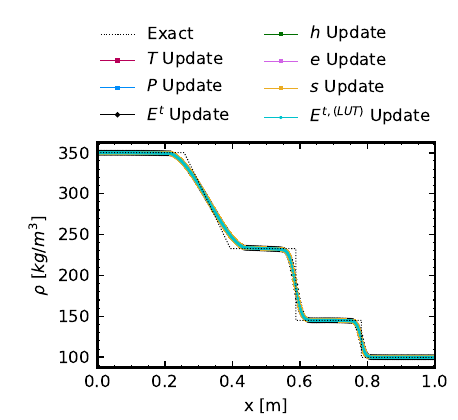}
        \label{subfig:co2_shock_comp_1st_rho}
    }%
    \subfloat[Pressure]{
        \includegraphics[scale=1]{
            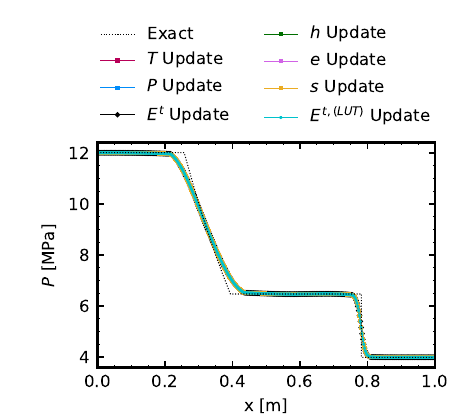}
        \label{subfig:co2_shock_comp_1st_p}
    } \\
    \subfloat[Density deviation from standard total energy update]{
        \includegraphics[scale=1]{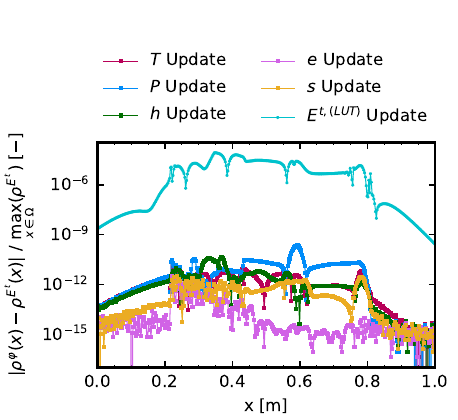}
        \label{subfig:co2_shock_err_1st_rho}
    }%
    \subfloat[Pressure deviation from standard total energy update]{
        \includegraphics[scale=1]{
            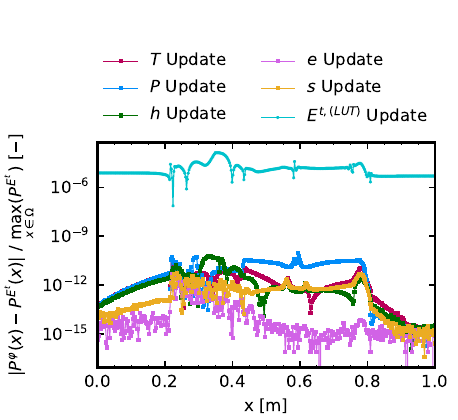}
        \label{subfig:co2_shock_err_1st_p}
    } \\
    \caption{Carbon-dioxide shock tube test 1st order results without using MUSCL with $\Delta x = 0.025\;\mathrm{m}$ using the Span-Wagner EoS~\cite{Span_1996} through the CoolProp thermodynamic library~\cite{CoolProp}. Density and pressure profiles (first row) and their dimensionless variations with respect to the standard total energy update (second row) at $t=0.001\;\mathrm{s}$. Comparison between the $\varphi$ update scheme for temperature, pressure, specific energy, specific enthalpy, specific entropy, and the standard total energy update (both using the analytical EoS and a look-up table).}
    \label{fig:co2_shock_comp_1st_order}
\end{figure}

In table~\ref{tab:secant_iterations} we report how many secant iterations per mesh element were needed to reach convergence in the $\varphibar$ search described in section~\ref{sec:varphibar_secant_search}. On average, the initial guess is already good enough to get close to machine precision, particularly in homogeneous portions of the domain like the post-shock plateau. Even in more demanding cases, such as the initial discontinuity, no more than eight iterations are needed for the $P$ update. Note that in this test we span thermodynamic states that are highly non-linear and therefore the maximum iteration counts shown here are most likely a worst case scenario.

\begin{table}[H]
    \centering
    \begin{tabular}{r |  c c}
            & \textbf{Average} & \textbf{Max} \\[1pt]
        \hline                                \\[-9pt]
        $T$ & $1.00473$        & $5$          \\[1pt]
        $P$ & $1.00451$        & $8$          \\[1pt]
        $e$ & $1.00376$        & $4$          \\[1pt]
        $h$ & $1.00263$        & $2$          \\[1pt]
        $s$ & $1.00279$        & $5$
    \end{tabular}
    \caption{Average and max (across all elements and time steps) number of secant iterations per element required to reach $\lvert F(\varphi) \rvert \le 10^{-14}$
        \label{tab:secant_iterations}}
\end{table}

\subsubsection{Computational costs}\label{sec:co2_computational_costs}

We now focus on analyzing the computational costs associated with the EoS evaluation. To do so we will categorize each EoS call performed during the simulation into three distinct types:
\begin{enumerate}
    \item \textbf{Riemann solver - HLLC}: computation of the pressure $P$, internal energy per unit mass $e$, speed of sound $c$ and the thermodynamic derivatives $\kappa = \left(\frac{\partial P}{\partial e}\right)_\rho$, $\chi = \left(\frac{\partial P}{\partial \rho}\right)_e$ at left and right of each face after MUSCL reconstruction.
    \item \textbf{Auxiliary}: computation of $P$ at left and right of each face after MUSCL reconstruction to check for unphysical reconstruction. Computation of $c$ on each element to compute the $\var{CFL}$ compliant time-step. Computation of $P$ on each element to check for unphysical solution updates.
    \item \textbf{Secant}: Computation of all thermodynamic quantities ($\Ephi$, $\Erho$, $E$) per element when performing the secant $\varphibar$ search. This cost does not exist by definition in the standard $\Et$ update.
\end{enumerate}

We then measure the computational time spent evaluating the EoS, and lump them in the aforementioned categories. To maintain the variability as low as possible we run these simulations serially on a single core at a locked frequency. Since some variability will inevitably appear, in Fig.~\ref{fig:co2_shock_computational_time} we present the average over $30$ runs and their respective ranges. All of the plots are scaled with respect to the standard $\Et$ results. In Fig.~\ref{subfig:co2_shock_time} we see the breakdown of the EoS evaluation times. The $T$ and $P$ update show significantly lower EoS evaluation times with speed-ups of up to $\sim 650\%$ for the $T$ update (see Fig.~\ref{subfig:co2_shock_time_total_speedup}), despite the additional $\sim 100\%$ EoS calls introduced by the secant search (see red bars in Fig.~\ref{subfig:co2_shock_time_eos_calls}). This is because, as visible in blue in Fig.~\ref{subfig:co2_shock_time_eos_calls}, each EoS call when using pressure or temperature is significantly faster than using the internal energy as done by the standard $\Et$ update. See appendix~\ref{appendix:n2_shock} for the same analysis in the dilute regime, where the speed-up for the $T$ update of around $700\%$ is comparable to the one measured here.

\begin{figure}
    \centering
    \subfloat[EoS evaluation time breakdown]{
        \includegraphics[scale=1]{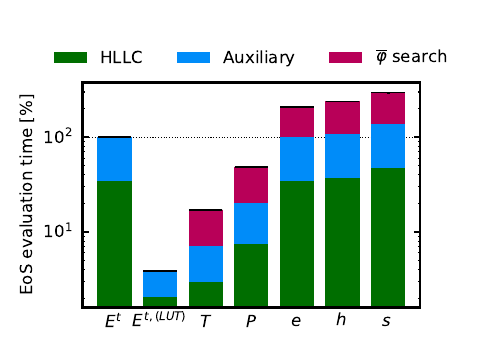}
        \label{subfig:co2_shock_time}
    }%
    \subfloat[EoS evaluation speed-up factor]{
        \includegraphics[scale=1]{
            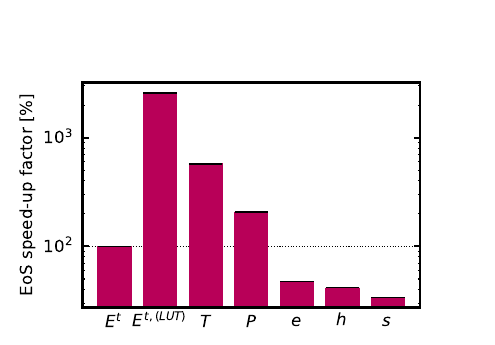}
        \label{subfig:co2_shock_time_speedup}
    }\\
    \subfloat[Total EoS calls and average EoS call time]{
        \includegraphics[scale=1]{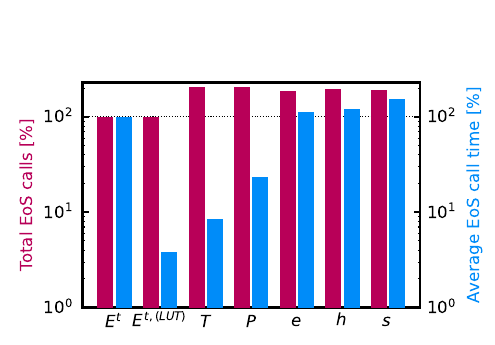}
        \label{subfig:co2_shock_time_eos_calls}
    }%
    \subfloat[Total speed-up factor]{
        \includegraphics[scale=1]{
            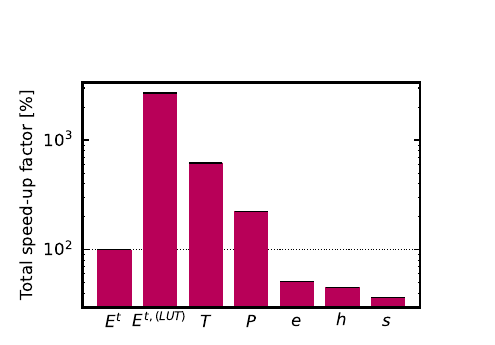}
        \label{subfig:co2_shock_time_total_speedup}
    }
    \caption{Carbon-dioxide shock tube test using MUSCL. EoS computational cost breakdown and comparison using the Span-Wagner EoS~\cite{Span_1996} through the CoolProp thermodynamic library~\cite{CoolProp}. All values are scaled with respect to the standard total energy update scheme, which is always $100\%$. The $y$ axis is in logarithmic scale. Comparison between the $\varphi$ update scheme for temperature, pressure, specific energy, specific enthalpy, specific entropy. We also show the computational costs obtained using CoolProp's tabular interpolation and the standard total energy update (both using the analytical EoS and a look-up table).}
    \label{fig:co2_shock_computational_time}
\end{figure}

\subsection{Nitrogen Shock Tube (VdW) - Shock speed and jump conditions}\label{sec:n2_shock}

In this section, we show the results for a nitrogen shock tube using the VdW EoS, see table~\ref{tab:n2_shock_data} for the data. This a standard Sod-type shock tube with a left-running rarefaction fan and right-running compression shock in dilute conditions. In appendix~\ref{appendix:n2_shock} we report the difference for different possible choices of $\varphi$ using the state-of-the-art thermodynamic library CoolProp~\cite{CoolProp}, which for nitrogen uses the Span-Wagner EoS~\cite{Span_2000}. Here we use the VdW EoS and the temperature $\varphi=T$ update to evaluate the shock speed and jump conditions in an over-refined simulation ($20000$ cells) and compare it to a simplified approximation and the exact analytical solution.

\begin{table}[H]
    \renewcommand{\arraystretch}{1.2}
    \centering
    \begin{tabular}{l l | l | l l}
        \textit{Equation of State}                    &                                     & \textit{Numerical}                   & ${x < 0\; \mathrm{m}}$         & ${x > 0\; \mathrm{m}}$         \\ \hline
        Van der Waals                                 &                                     & $x\in\left[-5,5\right]\; \mathrm{m}$ & $\rho=23.46 \;\mathrm{kg/m^3}$ & $\rho=11.73 \;\mathrm{kg/m^3}$ \\
        $a=173.943088\;\mathrm{m^5/kg\cdot s^2}$      & $\delta=0.4$                        & $t_{\var{final}}=0.01 \;\mathrm{s}$  & $u = 0 \; \mathrm{m/s}$        & $u = 0 \; \mathrm{m/s}$        \\
        $b=1.37851912 \cdot 10^{-3}\;\mathrm{m^3/kg}$ & $R=296.8\;\mathrm{{J}/{kg}\cdot K}$ & $\var{CFL}=0.1$                      & $P = 2 \; \mathrm{{MPa}}$      & $P = 1 \; \mathrm{{MPa}}$
    \end{tabular}
    \caption{Nitrogen shock tube test initial conditions, thermodynamic, domain and numerical data.}
    \label{tab:n2_shock_data}
\end{table}

We compare the $\varphibar$ search procedure presented in section~\ref{sec:varphibar_secant_search} to the approximate $\varphi$ update formula Eq~\ref{eq:varphi_update_formula_final} with $\rhobar=\rho^n$ and $\varphibar=\varphi^n$. In Fig.~\ref{fig:n2_shock_fixed_point} we see density, pressure, and velocity profiles compared to the exact solution, which at first glance appear identical. If we zoom on the right shock and on the velocity-pressure plateau in Fig.~\ref{fig:n2_shock_fixed_point_zoom} we understand why conservation is important. The results obtained using the approximate formula give wrong post-shock results (see Fig.~\ref{subfig:n2_shock_fixed_point_plateau_p}~\ref{subfig:n2_shock_fixed_point_plateau_u}) therefore jump conditions are not satisfied. Furthermore, the shock position is slightly to the right of the analytical solution (see Fig.~\ref{subfig:n2_shock_fixed_point_zoom_P}~\ref{subfig:n2_shock_fixed_point_zoom_u}). The spike we see in Fig.~\ref{subfig:n2_shock_fixed_point_plateau_p}~\ref{subfig:n2_shock_fixed_point_plateau_u} corresponds to the contact discontinuity, and is a feature of all Riemann solvers and is not caused by the time discretization. It is exaggerated by the tight zoom of Fig.~\ref{subfig:n2_shock_fixed_point_plateau_p}~\ref{subfig:n2_shock_fixed_point_plateau_u}. Looking at the absolute values of the total energy imbalance in time in Fig.~\ref{subfig:n2_shock_fixed_point_energy_conservation_log} it is clear how the approximate formula alone is not enough to conserve total energy, with a final value of roughly $10^9$ times greater than using the method proposed in this work.

\begin{figure}
    \centering
    \subfloat[Density]{
        \includegraphics[scale=1]{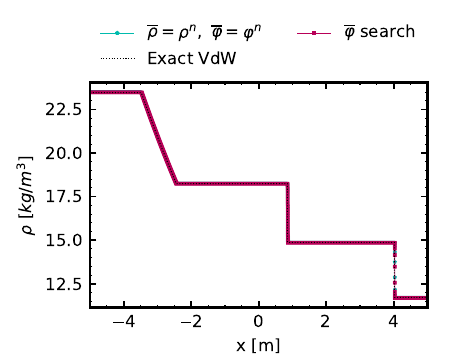}
        \label{subfig:n2_shock_fixed_point_rho}
    }%
    \subfloat[Pressure]{
        \includegraphics[scale=1]{
            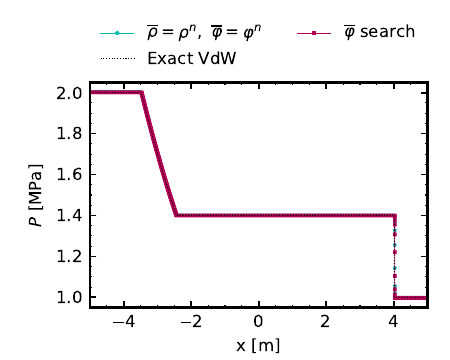}
        \label{subfig:n2_shock_fixed_point_p}
    }\\
    \subfloat[Velocity]{
        \includegraphics[scale=1]{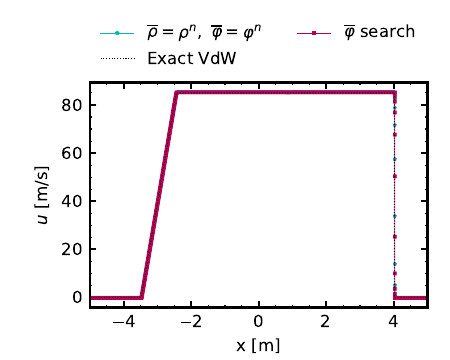}
        \label{subfig:n2_shock_fixed_point_u}
    }%
    \subfloat[Total energy imbalance]{
        \includegraphics[scale=1]{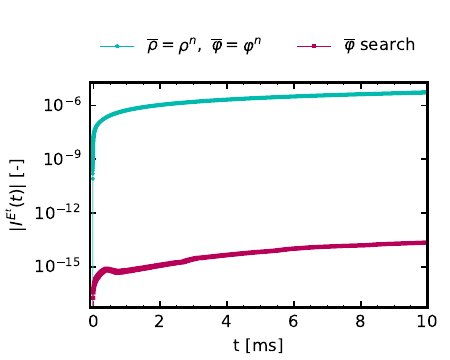}
        \label{subfig:n2_shock_fixed_point_energy_conservation_log}
    }
    \caption{Nitrogen shock tube test results using MUSCL with $\Delta x = 0.0002\;\mathrm{m}$ using the VdW EoS and the $\varphi=T$ update for temperature. Density, pressure, and velocity profiles at $t=0.01\;\mathrm{s}$ and absolute value of the total energy imbalance in time. Comparison between using $\rhobar=\left(\rho^n+\rho^{n+1}\right)/2$ with the secant search for $\varphibar=\Tbar$ and using a simplified approximation $\rhobar=\rho^n$ and $\varphibar=\varphi^n=T^n$ in the approximate $\varphi$ update formula.}
    \label{fig:n2_shock_fixed_point}
\end{figure}

\begin{figure}
    \centering
    \subfloat[Shock zoom - Pressure]{
        \includegraphics[scale=1]{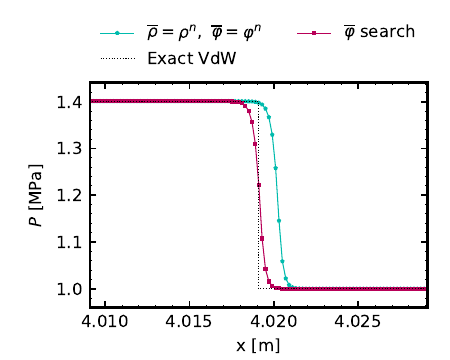}
        \label{subfig:n2_shock_fixed_point_zoom_P}
    }%
    \subfloat[Shock zoom - Velocity]{
        \includegraphics[scale=1]{
            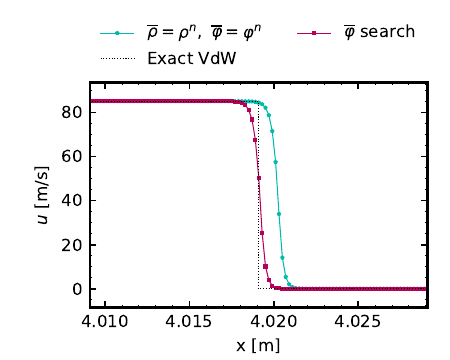}
        \label{subfig:n2_shock_fixed_point_zoom_u}
    }\\
    \subfloat[Plateau zoom - Pressure]{
        \includegraphics[scale=1]{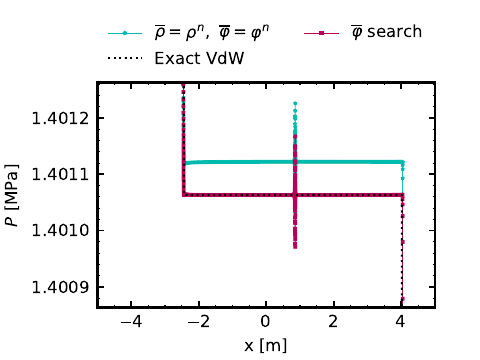}
        \label{subfig:n2_shock_fixed_point_plateau_p}
    }%
    \subfloat[Plateau zoom - Velocity]{
        \includegraphics[scale=1]{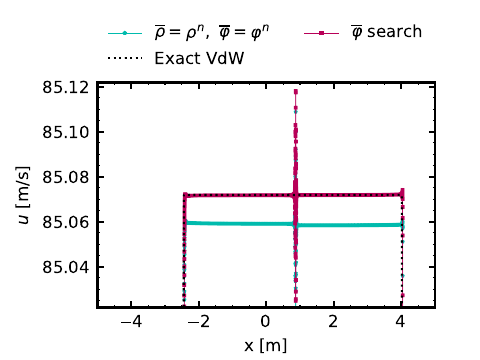}
        \label{subfig:n2_shock_fixed_point_plateau_u}
    }
    \caption{Nitrogen shock tube test results using MUSCL with $\Delta x = 0.0002\;\mathrm{m}$ using the VdW EoS and the $\varphi=T$ update for temperature. Pressure and velocity profiles at $t=0.01\;\mathrm{s}$.  Zoom on the right shock and on the plateaus for velocity and pressure. Comparison between using $\rhobar=\left(\rho^n+\rho^{n+1}\right)/2$ with the secant search for $\varphibar=\Tbar$ and using a simplified approximation $\rhobar=\rho^n$ and $\varphibar=\varphi^n=T^n$ in the approximate $\varphi$ update formula.}
    \label{fig:n2_shock_fixed_point_zoom}
\end{figure}

\subsection{Order Testing (VdW) - Method of Manufactured Solutions}\label{sec:mms}

We employ the method of manufactured solutions (MMS) to assess the spatial order of convergence. See Roache~\cite{Roache_1998} for an in-depth discussion on MMS. Using periodic boundary conditions in the $\Omega = [0,1]$ domain, we assume an exact periodic solution defined as
\begin{equation}
    \begin{dcases}
        \rho_{\textrm{MMS}}(x,t) & = \rho_0+ A_{\rho}\cos\left(\omega_{\rho,x}x+\omega_{\rho,t}t\right)     \\
        u_{\textrm{MMS}}(x,t)    & = u_0+ A_{u}\sin\left(\omega_{u,x}x+\omega_{u,t}t\right)                 \\
        T_{\textrm{MMS}}(x,t)    & = T_0+ A_{T}\sin\left(\omega_{T,x}x+\omega_{T,t}t\right)             \,.
    \end{dcases}
    \label{eq:exact_manufactured_solution}
\end{equation}

We substitute the manufactured solution Eq.~\eqref{eq:exact_manufactured_solution} into the Euler conservation equations Eq.~\eqref{eq:euler} to compute the source terms $\Theta$:
\begin{equation}
    \begin{dcases}
        \partiald{\rho_{\textrm{MMS}}(x,t)}{t}  & + \divergence \bigg[ \rhou_{\textrm{MMS}}(x,t) \bigg] = \Theta^\rho(x,t)                                                                     \\
        \partiald{\rhou_{\textrm{MMS}}(x,t)}{t} & + \divergence \bigg[ \rhou_{\textrm{MMS}}(x,t)\otimes \vec{u}_{\textrm{MMS}}(x,t) + P_{\textrm{MMS}}(x,t)\bigg] = \vec{\Theta}^{\rhou} (x,t) \\
        \partiald{\Et_{\textrm{MMS}}(x,t)}{t}   & + \divergence \bigg[ \left(\Et_{\textrm{MMS}}(x,t) + P_{\textrm{MMS}}(x,t)\right)\vec{u}_{\textrm{MMS}}(x,t) \bigg] =  \Theta^{\Et}(x,t)
    \end{dcases}
\end{equation}
using the following constants:
\begin{equation}
    \begin{aligned}
         & \rho_0 = 2\;\mathrm{kg/m^3} &  & \;\;\;\;A_{\rho} = 0.025\;\mathrm{kg/m^3} &  & \;\;\;\;\omega_{\rho,x} = 2\pi\;\mathrm{m^{-1}} &  & \;\;\;\;\omega_{\rho,t} = 256\pi\;\mathrm{s^{-1}}     \\
         & u_0 = 2\;\mathrm{m/s}       &  & \;\;\;\;A_{u} = 0.025\;\mathrm{m/s}       &  & \;\;\;\;\omega_{u,x} = 2\pi \;\mathrm{m^{-1}}   &  & \;\;\;\;\omega_{u,t} = 256\pi\;\mathrm{s^{-1}}        \\
         & T_0 = 300\;\mathrm{K}       &  & \;\;\;\;A_{T} = 0.1\;\mathrm{K}           &  & \;\;\;\;\omega_{T,x} = 2\pi\;\mathrm{m^{-1}}    &  & \;\;\;\;\omega_{T,t} = 256\pi\;\mathrm{s^{-1}}    \,.
    \end{aligned}
\end{equation}

We plot the dimensionless $L_2$ errors and the measured convergence rates of density, velocity, and temperature in Fig.~\ref{fig:mms1}, both with and without using MUSCL.  We also measure the error on the pressure. To adimensionalize we use the error of the coarsest mesh. As expected, first- and second-order convergence rates are measured, therefore the method proposed in this work preserves the expected characteristics of the underlying spatial discretization.

\begin{figure}
    \centering
    \subfloat[Error $L_2$ norm]{
        \includegraphics[scale=1]{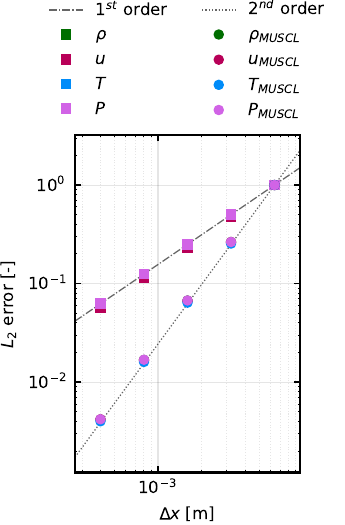}
        \label{subfig:mms1_l2_err}
    }%
    \subfloat[Measured convergence rate]{
        \includegraphics[scale=1]{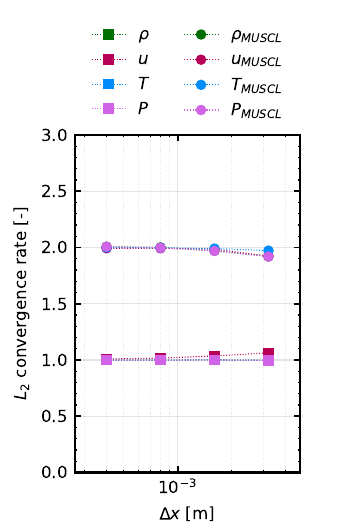}
        \label{subfig:mms1_l2_conv_rate}
    }
    \caption{MMS spatial convergence order test. Dimensionless $L_2$ norm of the errors and measured convergence rates at $t=0.0001s$. }
    \label{fig:mms1}
\end{figure}

\subsection{Dilute Gas Regime Tests (VdW)}\label{sec:classical}

In this section we show three tests in the dilute gas-dynamic regime using the VdW EoS. The first test is a 123 test (see table~\ref{tab:123_data}), whose solution is composed of two outgoing rarefaction waves. This test is useful to asses the behaviour of the scheme in the near vacuum region that forms as the rarefaction fans travel outwards. The second test is composed of two colliding shocks and is taken from Dumbser and Casulli~\cite{Dumbser_2016} (see table~\ref{tab:colliding_shock_data}). Here we can appreciate the behavior of the scheme in the presence of strong discontinuities and compressions. The third test is the so-called Lax shock which is also taken from~\cite{Dumbser_2016} (see table~\ref{tab:lax_shock_data}). The solution presents a large density jump in the central contact discontinuity with a strong right-running compression shock.  For all tests we will use the $T$ update and compare it to the standard $\Et$ update and the exact VdW solution.

\begin{table}[H]
    \renewcommand{\arraystretch}{1.2}
    \centering
    \begin{tabular}{l l | l | l l}
        \textit{Equation of State}         &                                   & \textit{Numerical}                  & ${x < 0.5\; \mathrm{m}}$   & ${x > 0.5\; \mathrm{m}}$   \\ \hline
        Van der Waals                      &                                   & $x\in\left[0,1\right]\; \mathrm{m}$ & $\rho=1 \;\mathrm{kg/m^3}$ & $\rho=1 \;\mathrm{kg/m^3}$ \\
        $a=0.5 \;\mathrm{m^5/kg\cdot s^2}$ & $\delta=0.4$                      & $t_{\var{final}}=0.2 \;\mathrm{s}$  & $u = -1 \; \mathrm{m/s}$   & $u = 1 \; \mathrm{m/s}$    \\
        $b=0.5 \;\mathrm{m^3/kg}$          & $R=0.4\;\mathrm{{J}/{kg}\cdot K}$ & $\var{CFL}=0.1$                     & $P = 0.4 \; \mathrm{{Pa}}$ & $P = 0.4 \; \mathrm{{Pa}}$
    \end{tabular}
    \caption{123 test initial conditions, thermodynamic, domain and numerical data.}
    \label{tab:123_data}
\end{table}

We report the density, pressure, and velocity profiles for the 123 test (see Fig.~\ref{fig:123_vdw}), colliding shocks (see Fig.~\ref{fig:colliding_shocks}) and Lax shock (see Fig.~\ref{fig:lax_shock}). Results compare well with the exact VdW solution everywhere, except for a kink in the density visible in Fig.~\ref{subfig:123_vdw_rho}, corresponding to the initial discontinuity in the 123 test. The $T$ update and the standard $\Et$ update show different behaviors here due to the fact that MUSCL is reconstructing and limiting two different variables. A similar behavior can be seen in the Lax shock on the pressure and velocity in Figs.~\ref{subfig:lax_shock_p}~\ref{subfig:lax_shock_u} across the contact discontinuity at $x\simeq0.65m$. The total energy imbalances in Figs.~\ref{subfig:123_vdw_energy_conservation}~\ref{subfig:colliding_shocks_energy_conservation}~\ref{subfig:lax_shock_energy_conservation} are almost overimposed for both the $T$ update and the standard $\Et$ update in all three test cases.

\begin{figure}
    \centering
    \subfloat[Density]{
        \includegraphics[scale=1]{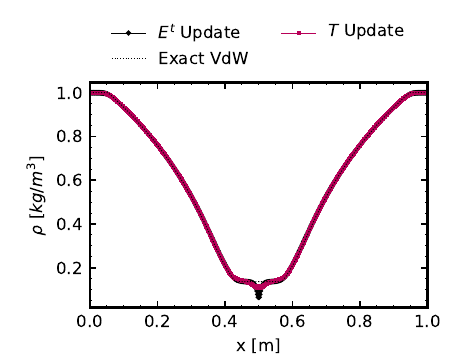}
        \label{subfig:123_vdw_rho}
    }%
    \subfloat[Pressure]{
        \includegraphics[scale=1]{
            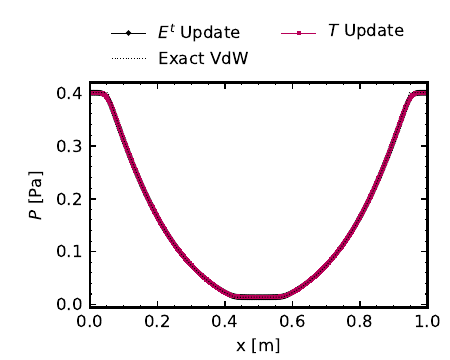}
        \label{subfig:123_vdw_p}
    }\\
    \subfloat[Velocity]{
        \includegraphics[scale=1]{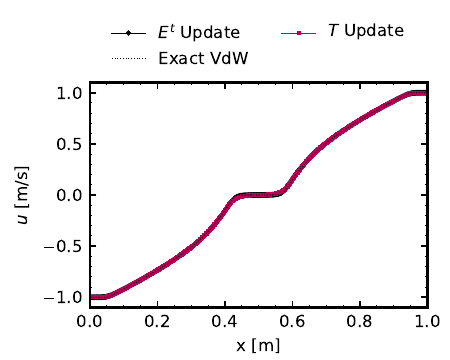}
        \label{subfig:123_vdw_u}
    }%
    \subfloat[Total energy imbalance]{
        \includegraphics[scale=1]{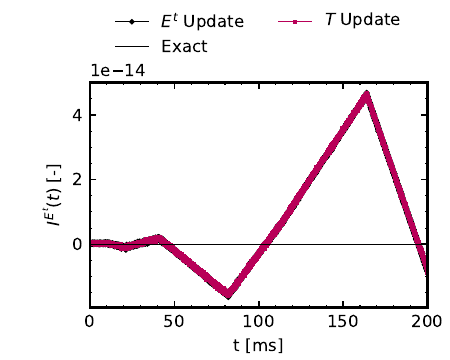}
        \label{subfig:123_vdw_energy_conservation}
    }
    \caption{123 test results using MUSCL with $\Delta x = 0.002\;\mathrm{m}$ using the VdW EoS. Density, pressure, and velocity profiles at $t=0.2\;\mathrm{s}$ and total energy imbalance in time. Comparison between the $\varphi$ update scheme for temperature and the standard total energy update.}
    \label{fig:123_vdw}
\end{figure}

\begin{table}[H]
    \renewcommand{\arraystretch}{1.2}
    \centering
    \begin{tabular}{l l | l | l l}
        \textit{Equation of State}         &                                   & \textit{Numerical}                  & ${x < 0.5\; \mathrm{m}}$   & ${x > 0.5\; \mathrm{m}}$   \\ \hline
        Van der Waals                      &                                   & $x\in\left[0,1\right]\; \mathrm{m}$ & $\rho=1 \;\mathrm{kg/m^3}$ & $\rho=1 \;\mathrm{kg/m^3}$ \\
        $a=0.5 \;\mathrm{m^5/kg\cdot s^2}$ & $\delta=0.4$                      & $t_{\var{final}}=0.1 \;\mathrm{s}$  & $u = 1 \; \mathrm{m/s}$    & $u = -1 \; \mathrm{m/s}$   \\
        $b=0.5 \;\mathrm{m^3/kg}$          & $R=0.4\;\mathrm{{J}/{kg}\cdot K}$ & $\var{CFL}=0.1$                     & $P = 2 \; \mathrm{{Pa}}$   & $P = 1 \; \mathrm{{Pa}}$
    \end{tabular}
    \caption{Colliding shocks test~\cite{Dumbser_2016} initial conditions, thermodynamic, domain and numerical data.}
    \label{tab:colliding_shock_data}
\end{table}

\begin{figure}
    \centering
    \subfloat[Density]{
        \includegraphics[scale=1]{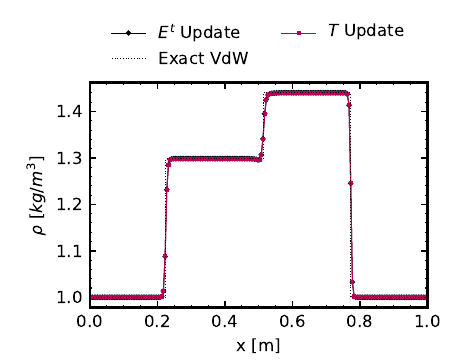}
        \label{subfig:colliding_shocks_rho}
    }%
    \subfloat[Pressure]{
        \includegraphics[scale=1]{
            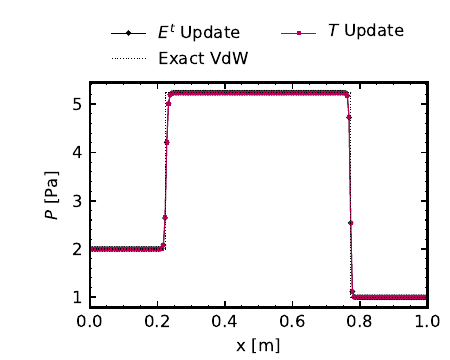}
        \label{subfig:colliding_shocks_p}
    }\\
    \subfloat[Velocity]{
        \includegraphics[scale=1]{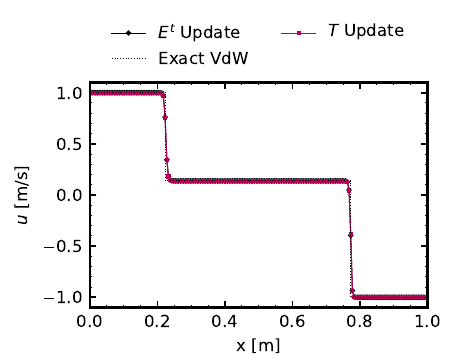}
        \label{subfig:colliding_shocks_u}
    }%
    \subfloat[Total energy imbalance]{
        \includegraphics[scale=1]{
            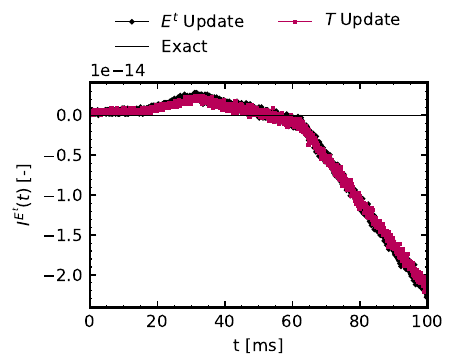}
        \label{subfig:colliding_shocks_energy_conservation}
    }
    \caption{Colliding shocks test~\cite{Dumbser_2016} results using MUSCL with $\Delta x = 0.005\;\mathrm{m}$ using the VdW EoS. Density, pressure, and velocity profiles at $t=0.1\;\mathrm{s}$ and total energy imbalance in time. Comparison between the $\varphi$ update scheme for temperature and the standard total energy update.}
    \label{fig:colliding_shocks}
\end{figure}

\begin{table}[H]
    \renewcommand{\arraystretch}{1.2}
    \centering
    \begin{tabular}{l l | l | l l}
        \textit{Equation of State}         &                                   & \textit{Numerical}                  & ${x < 0.5\; \mathrm{m}}$       & ${x > 0.5\; \mathrm{m}}$     \\ \hline
        Van der Waals                      &                                   & $x\in\left[0,1\right]\; \mathrm{m}$ & $\rho=0.445 \;\mathrm{kg/m^3}$ & $\rho=0.5 \;\mathrm{kg/m^3}$ \\
        $a=0.5 \;\mathrm{m^5/kg\cdot s^2}$ & $\delta=0.4$                      & $t_{\var{final}}=0.1 \;\mathrm{s}$  & $u = 0.698 \; \mathrm{m/s}$    & $u = 0 \; \mathrm{m/s}$      \\
        $b=0.5 \;\mathrm{m^3/kg}$          & $R=0.4\;\mathrm{{J}/{kg}\cdot K}$ & $\var{CFL}=0.1$                     & $P = 3.528 \; \mathrm{{Pa}}$   & $P = 0.571 \; \mathrm{{Pa}}$
    \end{tabular}
    \caption{Lax shock test initial conditions, thermodynamic, domain and numerical data.}
    \label{tab:lax_shock_data}
\end{table}

\begin{figure}
    \centering
    \subfloat[Density]{
        \includegraphics[scale=1]{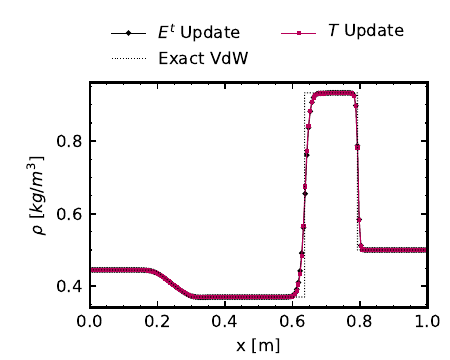}
        \label{subfig:lax_shock_rho}
    }%
    \subfloat[Pressure]{
        \includegraphics[scale=1]{
            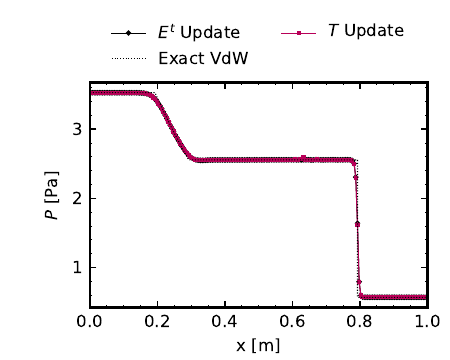}
        \label{subfig:lax_shock_p}
    }\\
    \subfloat[Velocity]{
        \includegraphics[scale=1]{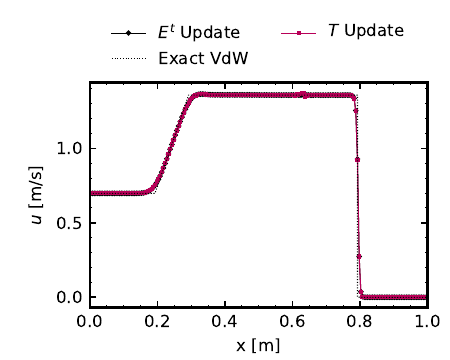}
        \label{subfig:lax_shock_u}
    }%
    \subfloat[Total energy imbalance]{
        \includegraphics[scale=1]{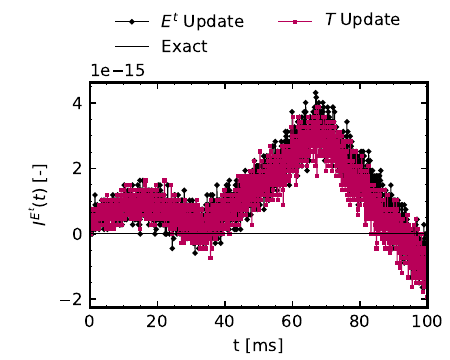}
        \label{subfig:lax_shock_energy_conservation}
    }
    \caption{Lax shock test results using MUSCL with $\Delta x = 0.005\;\mathrm{m}$ using the VdW EoS. Density, pressure and velocity profiles at $t=0.1\;\mathrm{s}$ and total energy imbalance in time. Comparison between the $\varphi$ update scheme for temperature and the standard total energy update.}
    \label{fig:lax_shock}
\end{figure}

\subsection{Dense Gas Regime Tests (VdW)}\label{sec:non_classical}

Here we show two dense gas tests named DG1 (see table~\ref{tab:dg1_data}) and DG2 (see table~\ref{tab:dg2_data}) from Guardone \& Vigevano~\cite{Guardone_2002}. These tests have been chosen to showcase the scheme in one of the most challenging gas-dynamic conditions, namely the non-classical regime~\cite{Bethe_1942,Zeldovich_1946,Thompson_1971} where rarefaction shock waves and compression fans are admissible solutions. These waves appear when the flow crosses a region of the thermodynamic plane where the fundamental derivative of gas-dynamics $\Gamma$, see Eq.~\eqref{eq:G_definition}, is less than $0$. This only happens for vapours of sufficient molecular complexity such as Siloxanes, see Colonna et al.~\cite{Colonna_2007}. When the flow crosses this region, the admissible solution switches from the classical \textit{rarefaction fan-compression shock} pair to the non-classical \textit{rarefaction shock-compression fan} pair. Furthermore, multiple crossing the $\Gamma=0$ line can lead to the rise of composite waves, see~\cite{Muller_2006,Zamfirescu_2008}. The first test DG1 is characterized by a wave structure sketched in Fig.~\ref{fig:dg1_wave_structure}, composed of a composite wave (rarefaction fan and non-classical rarefaction shock), a contact discontinuity, and a classical compression shock. The left and right states are in the convex region of the VdW EoS ($\Gamma>0$) but the fluid crosses the $\Gamma=0$ boundary during the evolution of the Riemann problem. The second test DG2 is characterized by a wave structure sketched in Fig.~\ref{fig:dg2_wave_structure}, composed by a non-classical rarefaction shock, a contact discontinuity, and a non-classical compression fan. For this test, the fluid always resides in the non-convex region of the VdW EoS ($\Gamma<0$)

\begin{table}[H]
    \renewcommand{\arraystretch}{1.2}
    \centering
    \begin{tabular}{l l | l | l l}
        \textit{Equation of State}   &                                              & \textit{Numerical}                  & ${x < 0.5\; \mathrm{m}}$      & ${x > 0.5\; \mathrm{m}}$       \\ \hline
        Van der Waals                & $\delta=0.0125$                              & $x\in\left[0,1\right]\; \mathrm{m}$ & $\rho=1.818\;\mathrm{kg/m^3}$ & $\rho=0.275 \;\mathrm{kg/m^3}$ \\
        $\rho_c=1 \;\mathrm{kg/m^3}$ & $P_c=1 \;\mathrm{Pa}$                        & $t_{\var{final}}=0.15 \;\mathrm{s}$ & $u = 0 \; \mathrm{m/s}$       & $u = 0 \; \mathrm{m/s}$        \\
        $T_c=1 \;\mathrm{K}$         & $R=2.\overline{6}\;\mathrm{{J}/{kg}\cdot K}$ & $\var{CFL}=0.1$                     & $P = 3 \; \mathrm{{Pa}}$      & $P = 0.575 \; \mathrm{{Pa}}$
    \end{tabular}
    \caption{DG1 dense gas test from~\cite{Guardone_2002} initial conditions, thermodynamic, domain and numerical data.}
    \label{tab:dg1_data}
\end{table}

\begin{figure}[H]
    \centering
    \includegraphics[scale=1]{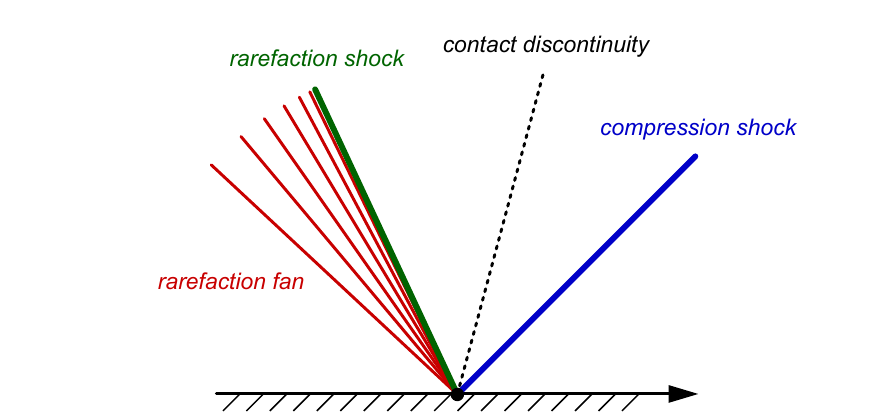}
    \caption{DG1 dense gas test from~\cite{Guardone_2002} solution wave structure sketch.}
    \label{fig:dg1_wave_structure}
\end{figure}

For both tests, we use the same $400$ element mesh as Guardone and Vigevano~\cite{Guardone_2002} and report the density, pressure, and fundamental derivative of gas-dynamics profiles for DG1 in Fig.~\ref{fig:dg1} and DG2 in Fig.~\ref{fig:dg2}. For the first test DG1, in Fig.~\ref{fig:dg1} we see that results are almost identical between the $T$ update, the standard $\Et$ update, and the results taken from Guardone and Vigevano~\cite{Guardone_2002} obtained using the method of Davis~\cite{Davis_1987}. The rarefaction shock at $x \simeq 0.53m$ is well caught with no spurious oscillations. For the second test DG2 in Fig.~\ref{fig:dg2}, both the rarefaction shock at $x \simeq 0.25m$ and the compression fan at $x \simeq 0.85m$ are well caught with no spurious oscillations. The results are in good agreement between the $T$ update, the standard $\Et$ update and the results taken from Guardone \& Vigevano~\cite{Guardone_2002} obtained using the method of Davis~\cite{Davis_1987}, except for a slight difference in the Davis results in the rightmost part of the compression fan at $x \simeq 0.88 m$. The total energy imbalance for both DG1 and DG2 tests are almost overimposed between the $T$ update and the standard $\Et$ update as seen in Figs.~\ref{subfig:dg1_energy_conservation}~\ref{subfig:dg2_energy_conservation}. The fundamental derivative of gas-dynamics for test DG2 in Fig.~\ref{subfig:dg2_G} is negative everywhere, while for test DG1 in Fig.~\ref{subfig:dg1_G} it only crosses the $\Gamma=0$ boundary during the simulation, remaining positive everywhere else.

\begin{figure}
    \centering
    \subfloat[Density]{
        \includegraphics[scale=1]{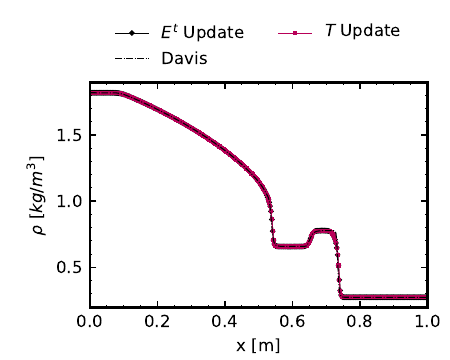}
        \label{subfig:dg1_rho}
    }%
    \subfloat[Pressure]{
        \includegraphics[scale=1]{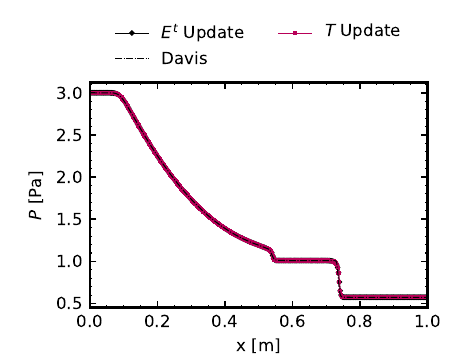}
        \label{subfig:dg1_p}
    }\\
    \subfloat[Fundamental derivative of gas-dynamics]{
        \includegraphics[scale=1]{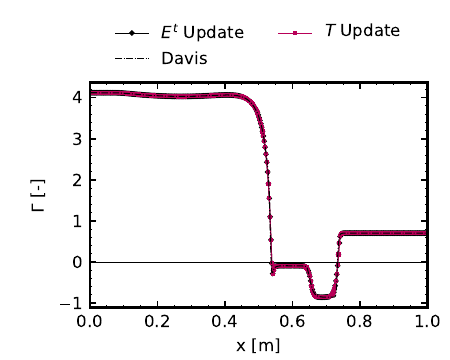}
        \label{subfig:dg1_G}
    }%
    \subfloat[Total energy imbalance]{
        \includegraphics[scale=1]{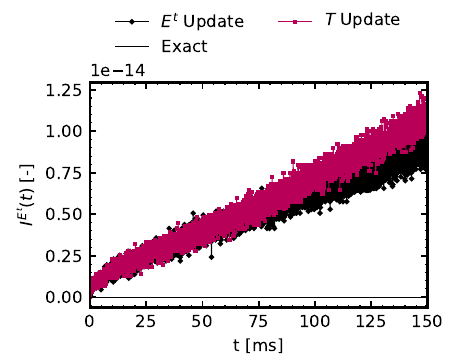}
        \label{subfig:dg1_energy_conservation}
    }
    \caption{DG1 dense gas~\cite{Guardone_2002} results using MUSCL with $\Delta x = 0.0025\;\mathrm{m}$ using the VdW EoS. Density, pressure and fundamental derivative of gas-dynamics profiles at $t=0.15\;\mathrm{s}$ and total energy imbalance in time. Comparison between the $\varphi$ update scheme for temperature and the standard total energy update. Results with Davis~\cite{Davis_1987} method taken from Guardone \& Vigevano~\cite{Guardone_2002}.}
    \label{fig:dg1}
\end{figure}

\begin{table}[H]
    \renewcommand{\arraystretch}{1.2}
    \centering
    \begin{tabular}{l l | l | l l}
        \textit{Equation of State}   &                                              & \textit{Numerical}                  & ${x < 0.5\; \mathrm{m}}$      & ${x > 0.5\; \mathrm{m}}$       \\ \hline
        Van der Waals                & $\delta=0.0125$                              & $x\in\left[0,1\right]\; \mathrm{m}$ & $\rho=0.879\;\mathrm{kg/m^3}$ & $\rho=0.562 \;\mathrm{kg/m^3}$ \\
        $\rho_c=1 \;\mathrm{kg/m^3}$ & $P_c=1 \;\mathrm{Pa}$                        & $t_{\var{final}}=0.45 \;\mathrm{s}$ & $u = 0 \; \mathrm{m/s}$       & $u = 0 \; \mathrm{m/s}$        \\
        $T_c=1 \;\mathrm{K}$         & $R=2.\overline{6}\;\mathrm{{J}/{kg}\cdot K}$ & $\var{CFL}=0.1$                     & $P = 1.09 \; \mathrm{{Pa}}$   & $P = 0.885 \; \mathrm{{Pa}}$
    \end{tabular}
    \caption{DG2 dense gas test from~\cite{Guardone_2002} initial conditions, thermodynamic, domain and numerical data.}
    \label{tab:dg2_data}
\end{table}

\begin{figure}
    \centering
    \includegraphics[scale=1]{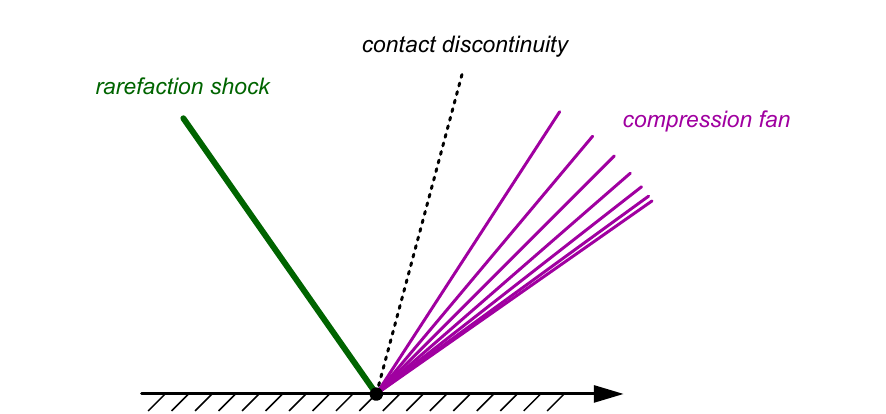}
    \caption{DG2 dense gas test from~\cite{Guardone_2002} solution wave structure.}
    \label{fig:dg2_wave_structure}
\end{figure}

\begin{figure}
    \centering
    \subfloat[Density]{
        \includegraphics[scale=1]{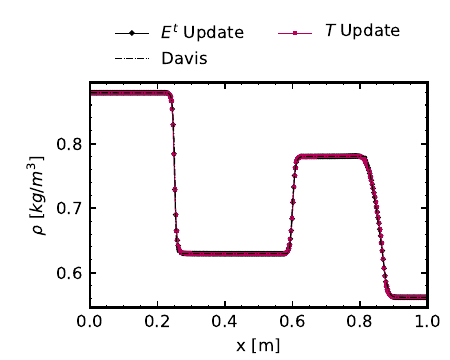}
        \label{subfig:dg2_rho}
    }%
    \subfloat[Pressure]{
        \includegraphics[scale=1]{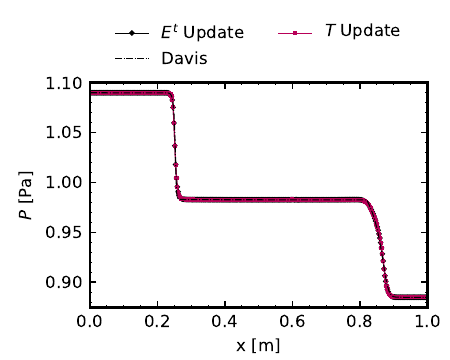}
        \label{subfig:dg2_p}
    }\\
    \subfloat[Fundamental derivative of gas-dynamics]{
        \includegraphics[scale=1]{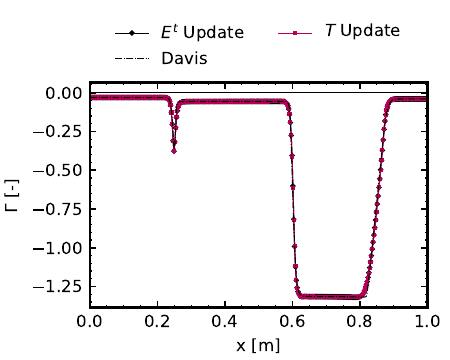}
        \label{subfig:dg2_G}
    }%
    \subfloat[Total energy imbalance]{
        \includegraphics[scale=1]{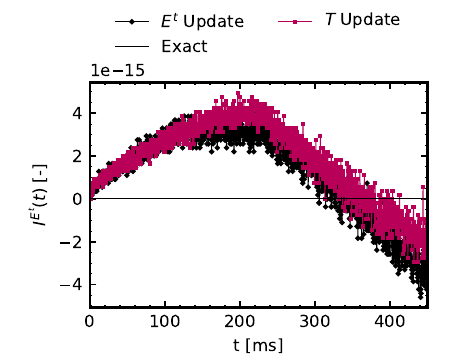}
        \label{subfig:dg2_energy_conservation}
    }
    \caption{DG2 non-ideal~\cite{Guardone_2002} results using MUSCL with $\Delta x = 0.0025\;\mathrm{m}$ using the VdW EoS. Density, pressure and fundamental derivative of gas-dynamics profiles at $t=0.45\;\mathrm{s}$ and total energy imbalance in time. Comparison between the $\varphi$ update scheme for temperature and the standard total energy update. Results with Davis~\cite{Davis_1987} method taken from Guardone \& Vigevano~\cite{Guardone_2002}.}
    \label{fig:dg2}
\end{figure}

\section{Conclusions}\label{sec:conclusions}

We presented a procedure to solve the Euler equations by updating any thermodynamic variable instead of the total energy in a conservative manner. The procedure is agnostic to the chosen equation of state and spatial discretization. We measured an increase in the number of total thermodynamic evaluations of $70-100\%$ when compared to the standard total energy update. The computational cost of this increase can be completely offset if using the temperature as the new thermodynamic variable because each thermodynamic evaluation is significantly faster when using the analytical Span-Wagner EoS. We measured an average speed-up of $650-700\%$ if compared to the standard total energy update. Results obtained by updating temperature, pressure, specific energy, specific enthalpy, or specific entropy instead of the total energy are all within the margin of error, with all of them conserving total energy close to machine precision and to the standard total energy update. We showed that using a limited MUSCL slope reconstruction the measured order of spatial convergence remains $\Delta x^2$ as expected. Although we only presented 1D results, the procedure is applicable to 2D or 3D since the increase in dimensions would not affect the time update in any way as it only looks at one element. The thermodynamic evaluation speed-up should be even larger in 3D as the main number of thermodynamic evaluations happens at the faces, the number of which scales faster than the number of control volumes. Adding viscous or turbulent diffusion effects to the spatial residuals does not change the derivation and application of the method. The scheme in its current form is limited to first order in time and is still CFL bound as the standard forward Euler approach. Higher order spatial discretizations such as WENO reconstruction would not affect the time update scheme presented, while approaches such as discontinuous-Galerkin could require some additional attention, in particular to the spatial integrals appearing in the approximate update formula. On the other hand, higher order time discretizations would require non-trivial modifications depending on the chosen scheme. The iterative procedure that retrieves conservation is limited to using a secant root search algorithm as the derivatives of the function $F(\varphibar)$ are not easily computable. Studying the application of more efficient root search algorithms could improve the performance of the scheme. Future work will also focus on the extension to multi-phase simulations.

\section*{Acknowledgements}
B. Re acknowledges the financial support received by the European Union's Horizon Europe programme under the MSCA-2021-PF grant agreement Id 101066019 (project NI2PhORC). G. Sirianni acknowledges the financial support of the European Union under the NextGenerationEU action.

\section*{Dataset}
A dataset containing all the numerical results described in this paper is available on Zenodo~\cite{zenodo_dataset}.

\bibliographystyle{unsrt}

\bibliography{bibliography}

\newpage

\appendix

\section{Polytropic Van der Waals Equation of State}\label{appendix:vdw}

We report here useful formulas and thermodynamic derivatives of the polytropic VdW EoS for the $T$ update. Pressure, internal energy, and the fundamental derivative of gas dynamics are:

\begin{equation}
    \begin{aligned}
        P(\rho, T)     & =  \dfrac{\rho R T}{1-b\rho} - a\rho^2                                                                                                          \\[10pt]
        E(\rho, T)     & = \dfrac{\rho R T}{\delta}  - a\rho^2                                                                                                           \\
        \Gamma(\rho,P) & = \dfrac{(\delta+1)(\delta+2) \dfrac{P+a \rho^2}{(1-\rho b)^2}\rho^2-6 a \rho^4}{2(\delta+1)\dfrac{P+a \rho^2}{(1-\rho b) }\rho^2 - 4 a \rho^4}
    \end{aligned}
\end{equation}

Other miscellaneous quantities are:
\begin{equation}
    \begin{aligned}[c]
         & a = \dfrac{27}{64}\dfrac{\left(R T_c \right)^2}{P_c} &  & b = \dfrac{1}{8}\dfrac{ R T_c  }{P_c}                                                          \\[10pt]
         & \kappa = \dfrac{\delta}{1-b \rho}                    &  & \chi = \delta \dfrac{b\left(E - a\rho^2\right) + 2 a \rho}{\left(1-b \rho\right)^2} - 2 a \rho
    \end{aligned}
    \label{eq:vdw}
\end{equation}

The internal energy derivatives needed for the $\varphi$ update scheme for temperature $T$:
\begin{equation}
    \begin{aligned}
         & \Ephi = \left(\partiald{E}{T}\right)_\rho = \dfrac{\rho R}{\delta}             \\[10pt]
         & \Erho = \left(\partiald{E}{\rho}\right)_T = \dfrac{R T }{\delta } - 2 a   \rho
    \end{aligned}
\end{equation}

\section{Nitrogen Shock Tube with Span-Wagner EoS}\label{appendix:n2_shock}

Results obtained using the CoolProp thermodynamic library, leveraging the Span-Wagner EoS~\cite{Span_2000}. In Fig.~\ref{fig:n2_shock_comp} we show profiles of density, pressure, and velocity, compared to the exact solution for nitrogen computed using the VdW EoS. There is virtually no difference in the results between all $\varphi$ choices, and they all overlap perfectly with the solution obtained using the standard $\Et$ update. If we focus on the total energy imbalance in Fig.~\ref{subfig:n2_shock_comp_energy_conservation_log} we see how some choices of $\varphi$ behave slightly differently from others. Namely, the specific internal energy $e$ behaves exactly identically to the standard $\Et$ update, while temperature $T$ and specific enthalpy $h$ show some similarity. Pressure $P$ and specific entropy $s$ showcase some sporadic steps that may at first glance appear large, but are still of the order of $10^{-13}$.

\begin{figure}
    \centering
    \subfloat[Density]{
        \includegraphics[scale=1]{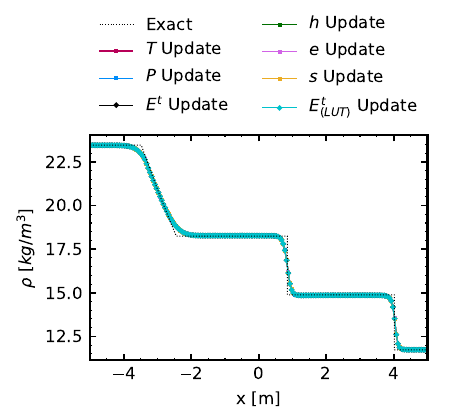}
        \label{subfig:n2_shock_comp_rho}
    }%
    \subfloat[Pressure]{
        \includegraphics[scale=1]{
            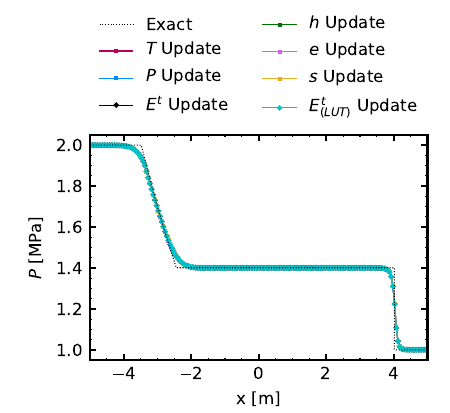}
        \label{subfig:n2_shock_comp_p}
    }\\
    \subfloat[Velocity]{
        \includegraphics[scale=1]{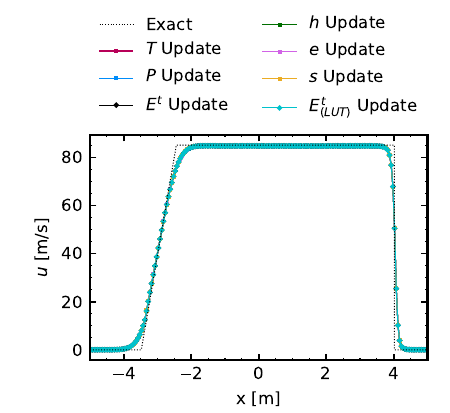}
        \label{subfig:n2_shock_comp_u}
    }%
    \subfloat[Total energy imbalance]{
        \includegraphics[scale=1]{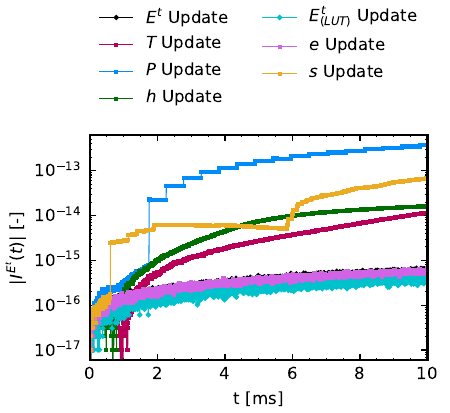}
        \label{subfig:n2_shock_comp_energy_conservation_log}
    }
    \caption{Nitrogen shock tube test results with $\Delta x = 0.05\;\mathrm{m}$ using the Span-Wagner EoS~\cite{Span_2000} through the CoolProp thermodynamic library~\cite{CoolProp}. Density, pressure and velocity profiles at $t=0.01\;\mathrm{s}$ and absolute value of the total energy imbalance in time. Comparison between the $\varphi$ update scheme for temperature, pressure, specific energy, specific enthalpy, specific entropy and the standard total energy update.}
    \label{fig:n2_shock_comp}
\end{figure}

\begin{figure}
    \centering
    \subfloat[EoS evaluation time breakdown]{
        \includegraphics[scale=1]{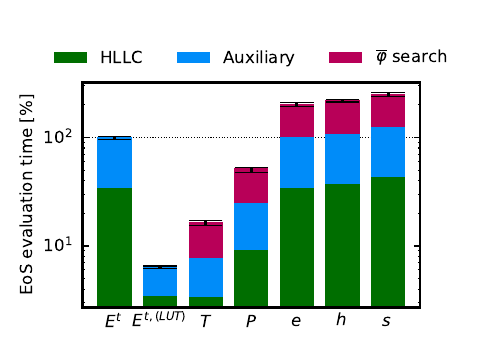}
        \label{subfig:n2_shock_time}
    }%
    \subfloat[EoS evaluation speed-up factor]{
        \includegraphics[scale=1]{
            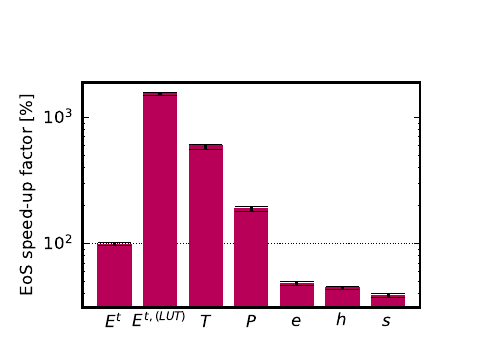}
        \label{subfig:n2_shock_time_speedup}
    }\\
    \subfloat[Total EoS calls and average EoS call time]{
        \includegraphics[scale=1]{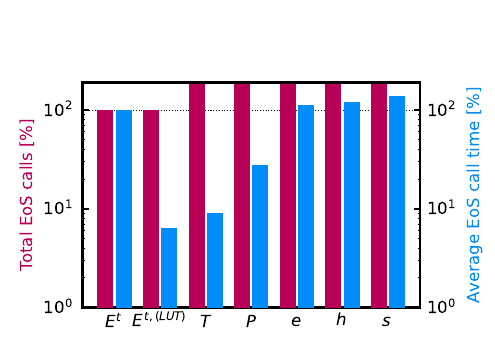}
        \label{subfig:n2_shock_time_eos_calls}
    }%
    \subfloat[Total speed-up factor]{
        \includegraphics[scale=1]{
            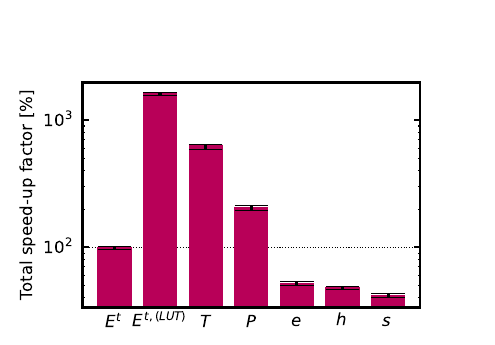}
        \label{subfig:n2_shock_time_total_speedup}
    }
    \caption{Nitrogen shock tube test. EoS computational cost breakdown and comparison using the Span-Wagner EoS~\cite{Span_2000} through the CoolProp thermodynamic library~\cite{CoolProp}. All values are scaled with respect to the standard total energy update scheme, which is always $100\%$. The $y$ axis is in logarithmic scale. Comparison between the $\varphi$ update scheme for temperature, pressure, specific energy, specific enthalpy, and specific entropy. We also show the computational costs obtained using CoolProp's tabular interpolation and the standard total energy update.}
    \label{fig:n2_shock_computational_time}
\end{figure}

\section{Nitrogen Shock Tube (VdW) - Choice of Linearization Density}\label{sec:n2_shock_lin_rho}

We show here how the choice of linearization density does not affect the accuracy of the solution. Results for the nitrogen shock tube with the VdW EoS from section~\ref{sec:n2_shock} are reported here using the $\varphibar$ search, but with three different choices for $\rhobar$, namely $\rho^n$, $\rho^{n+1}$ and $(\rho^{n} + \rho^{n+1})/2$. The results are reported in Fig.~\ref{fig:n2_shock_fixed_point_zoom_rhon_rhon1}.

\begin{figure}
    \centering
    \subfloat[Shock zoom - Pressure]{
        \includegraphics[scale=1]{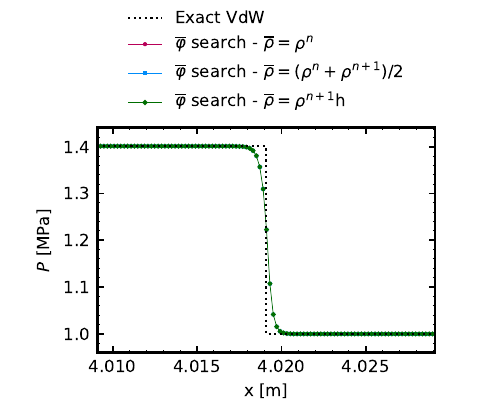}
        \label{subfig:n2_shock_fixed_point_rhon_rhon1_zoom_P}
    }%
    \subfloat[Shock zoom - Velocity]{
        \includegraphics[scale=1]{
            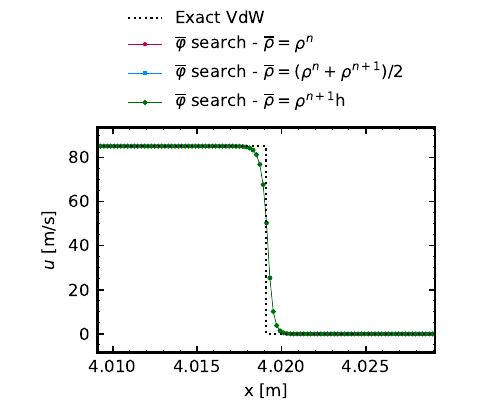}
        \label{subfig:n2_shock_fixed_point_rhon_rhon1_zoom_u}
    }\\
    \subfloat[Plateau zoom - Pressure]{
        \includegraphics[scale=1]{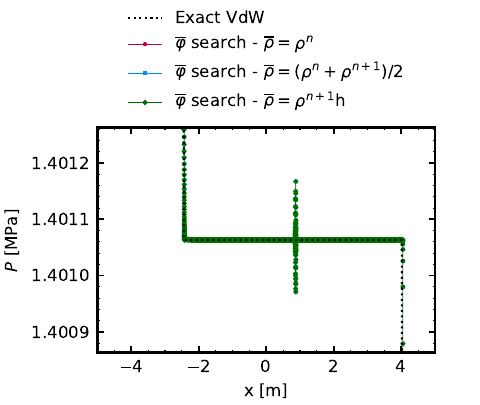}
        \label{subfig:n2_shock_fixed_point_rhon_rhon1_plateau_p}
    }%
    \subfloat[Plateau zoom - Velocity]{
        \includegraphics[scale=1]{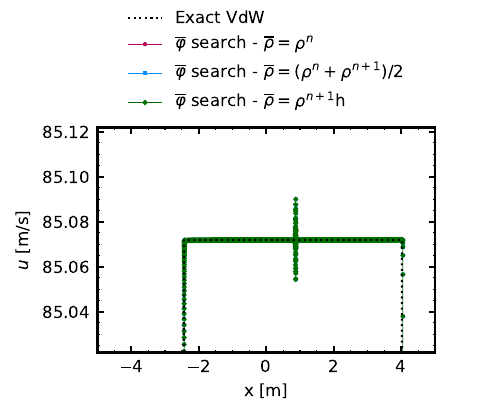}
        \label{subfig:n2_shock_fixed_point_rhon_rhon1_plateau_u}
    }
    \caption{Nitrogen shock tube test results using MUSCL with $\Delta x = 0.0002\;\mathrm{m}$ using the VdW EoS and the $\varphi=T$ update for temperature. Pressure and velocity profiles at $t=0.01\;\mathrm{s}$. Zoom on the right shock and on the plateaus for velocity and pressure. Comparison between using $\rhobar=\left(\rho^n+\rho^{n+1}\right)/2$, $\rhobar=\rho^n$, and $\rhobar=\rho^{n+1}$ with the secant search for $\varphibar=\Tbar$.}
    \label{fig:n2_shock_fixed_point_zoom_rhon_rhon1}
\end{figure}

\end{document}